\documentclass[11pt,leqno]{article}
\usepackage{amsfonts}
\usepackage{graphicx}

\begin{document}

\def\R{{\mathbb R}}
\def\T{{\mathbb T}}
\def\S{{\mathbb S}}
\def\C{{\mathbb C}}
\def\Z{{\mathbb Z}}
\def\N{{\mathbb N}}
\def\H{{\mathbb H}}
\def\B{{\mathbb B}}
\def\diam{\mbox{\rm diam}}
\def\sn{\S^{n-1}}
\def\rr{{\cal R}}
\def\mt{{\Lambda}}
\def\e{\emptyset}
\def\dQ{\partial Q}
\def\dk{\partial K}
\def\endofproof{{\rule{6pt}{6pt}}}
\def\di{\displaystyle}
\def\dist{\mbox{\rm dist}}
\def\sa+{\Sigma_A^+}
\def\du{\frac{\partial}{\partial u}}
\def\dv{\frac{\partial}{\partial v}}
\def\dt{\frac{d}{d t}}
\def\dx{\frac{\partial}{\partial x}}
\def\con{\mbox{\rm const }}
\def\nn{{\cal N}}
\def\mm{{\cal M}}
\def\kk{{\cal K}}
\def\ll{{\cal L}}
\def\vv{{\cal V}}
\def\bb{{\cal B}}
\def\ma{\mm_{a}}
\def\lab{L_{ab}}
\def\mabn{\mm_{a}^N}
\def\man{\mm_a^N}
\def\labn{L_{ab}^N}
\def\fa{f^{(a)}}
\def\ff{{\cal F}}
\def\i{{\bf i}}
\def\gge{{\cal G}_\epsilon}
\def\gej{\chi^{(j)}_\mu}
\def\ge{\chi_\epsilon}
\def\geo{\chi^{(1)}_\mu}
\def\get{\chi^{(2)}_\mu}
\def\gei{\chi^{(i)}_{\mu}}
\def\gee{\chi_{\mu}}
\def\gett{\chi^{(2)}_{\mu}}
\def\geol{\chi^{(1)}_{\ell}}
\def\getl{\chi^{(2)}_{\ell}}
\def\geil{\chi^{(i)}_{\ell}}
\def\gee{\chi_{\ell}}
\def\tt{{\cal T}}
\def\uu{{\cal U}}
\def\wloc{W_{\epsilon}}
\def\Int{\mbox{\rm Int}}
\def\dist{\mbox{\rm dist}}
\def\pr{\mbox{\rm pr}}
\def\pp{{\cal P}}
\def\aa{{\cal A}}
\def\cc{{\cal C}}
\def\supp{\mbox{\rm supp}}
\def\Arg{\mbox{\rm Arg}}
\def\In{\mbox{\rm Int}}
\def\con{\mbox{\rm const}\;}
\def\Re{\mbox{\rm Re}}
\def\li{\mbox{\rm li}} 
\def\Seo{S^*_\epsilon(\Omega)}
\def\sdk{S^*_{\dk}(\Omega)}
\def\lae{\Lambda_{\epsilon}}
\def\ep{\epsilon}
\def\oo{{\cal O}}
\def\be{\begin{equation}}
\def\ee{\end{equation}}
\def\beqn{\begin{eqnarray*}}
\def\eeqn{\end{eqnarray*}}
\def\Pr{\mbox{\rm Pr}}

\def\gi{\gamma^{(i)}}
\def\ii{{\imath }}
\def\jj{{\jmath }}
\def\II{{\cal I}}
\def\ccij{ \cc_{i'_0,j'_0}[\eta]}
\def\dd{{\cal D}}
\def\la{\langle}
\def\ra{\rangle}
\def\bs{\bigskip}
\def\xio{\xi^{(0)}}
\def\xo{x^{(0)}}
\def\zo{z^{(0)}}
\def\Con{\mbox{\rm Const}\;}
\def\do{\partial \Omega}
\def\dk{\partial K}
\def\dl{\partial L}
\def\ll{{\cal L}}
\def\kk{{\cal K}}
\def\kk{{\cal K}}
\def\pr{{\rm pr}}
\def\ff{{\cal F}}
\def\G{{\cal G}}
\def\C{{\bf C}}
\def\dist{{\rm dist}}
\def\dds{\frac{d}{ds}}
\def\con{{\rm const}\;}
\def\Con{{\rm Const}\;}
\def\di{\displaystyle}
\def\oo{\mbox{\rm O}}
\def\hess{\mbox{\rm Hess}}
\def\gi{\gamma^{(i)}}
\def\endofproof{{\rule{6pt}{6pt}}}
\def\xm{x^{(m)}}
\def\vm{\varphi^{(m)}}
\def\km{k^{(m)}}
\def\dm{d^{(m)}}
\def\kam{\kappa^{(m)}}
\def\dem{\delta^{(m)}}
\def\xim{\xi^{(m)}}
\def\ep{\epsilon}
\def\ms{\medskip}
\def\ex{\mbox{\rm extd}}

\def\clip{C^{\mbox{\footnotesize \rm Lip}}}
\def\wlocs{W^s_{\mbox{\footnote\rm loc}}}
\def\Lip{\mbox{\rm Lip}}

\def\Xr{X^{(r)}}
\def\lip{\mbox{{\footnotesize\rm Lip}}}
\def\Vol{\mbox{\rm Vol}}

\def\naf{\nabla f(z)}
\def\so{\sigma_0}
\def\Xo{X^{(0)}}
\def\z1{z^{(1)}}
\def\Vo{V^{(0)}}
\def\Yo{Y{(0)}}

\def\uo{u^{(0)}}
\def\vo{v^{(0)}}
\def\no{\nu^{(0)}}
\def\psa{\partial^{(s)}_a}
\def\hcd{\hc^{(\delta)}}
\def\Md{M^{(\delta)}}
\def\Uo{U^{(1)}}
\def\Ut{U^{(2)}}
\def\Uj{U^{(j)}}
\def\no{n^{(1)}}
\def\nt{n^{(2)}}
\def\nj{n^{(j)}}
\def\ccm{\cc^{(m)}}

\def\ooo{\oo^{(1)}}
\def\oot{\oo^{(2)}}
\def\ooj{\oo^{(j)}}
\def\fo{f^{(1)}}
\def\ft{f^{(2)}}
\def\fj{f^{(j)}}
\def\wo{w^{(1)}}
\def\wt{w^{(2)}}
\def\wj{w^{(j)}}
\def\Vo{V^{(1)}}
\def\Vt{V^{(2)}}
\def\Vj{V^{(j)}}

\def\Ul{U^{(\ell)}}
\def\Uj{U^{(j)}}
\def\wl{w^{(\ell)}}
\def\Vl{V^{(\ell)}}
\def\Ujj{U^{(j+1)}}
\def\wjj{w^{(j+1)}}
\def\Vjj{V^{(j+1)}}
\def\Ujo{U^{(j_0)}}
\def\wjo{w^{(j_0)}}
\def\Vjo{V^{(j_0)}}
\def\vj{v^{(j)}}
\def\vl{v^{(\ell)}}

\def\f0{f^{(0)}}

\def\gl{\gamma_\ell}
\def\id{\mbox{\rm id}}
\def\piU{\pi^{(U)}}

\def\cca{C^{(a)}}
\def\bba{B^{(a)}}
\def\co{\; \stackrel{\circ}{C}}

\def\oV{\overline{V}}
\def\saa{\Sigma^+_A}
\def\sa{\Sigma_A}
\def\mta{\Lambda(A, \tau)}
\def\mtaa{\Lambda^+(A, \tau)}

\def\Int{\mbox{\rm Int}}
\def\epo{\ep^{(0)}}
\def\pH{\partial \H^{n+1}}
\def\sh{S^*(\H^{n+1})}
\def\zoo{z^{(1)}}
\def\yoo{y^{(1)}}
\def\xoo{x^{(1)}}


\def\supp{\mbox{\rm supp}}
\def\Arg{\mbox{\rm Arg}}
\def\In{\mbox{\rm Int}}
\def\diam{\mbox{\rm diam}}
\def\e{\emptyset}
\def\endofproof{{\rule{6pt}{6pt}}}
\def\di{\displaystyle}
\def\dist{\mbox{\rm dist}}
\def\con{\mbox{\rm const }}
\def\Box{\spadesuit}
\def\Int{\mbox{\rm Int}}
\def\dist{\mbox{\rm dist}}
\def\pr{\mbox{\rm pr}}
\def\be{\begin{equation}}
\def\ee{\end{equation}}
\def\beqn{\begin{eqnarray*}}
\def\eeqn{\end{eqnarray*}}
\def\la{\langle}
\def\ra{\rangle}
\def\bs{\bigskip}
\def\Con{\mbox{\rm Const}\;}
\def\clip{C^{\mbox{\footnotesize \rm Lip}}}
\def\wlocs{W^s_{\mbox{\footnote\rm loc}}}
\def\Lip{\mbox{\rm Lip}}
\def\lip{\mbox{\footnotesize\rm Lip}}
\def\Re{\mbox{\rm Re}}
\def\li{\mbox{\rm li}} 
\def\ep{\epsilon}
\def\ms{\medskip}
\def\dds{\frac{d}{ds}}
\def\oo{\mbox{\rm O}}
\def\hess{\mbox{\rm Hess}}
\def\id{\mbox{\rm id}}
\def\ii{{\imath }}
\def\jj{{\jmath }}
\def\graph{\mbox{\rm graph}}
\def\span{\mbox{\rm span}}

\def\i{{\bf i}}
\def\C{{\bf C}}

\def\ss{{\cal S}}
\def\tt{{\cal T}}
\def\E{{\cal E}}
\def\rr{{\cal R}}
\def\nn{{\cal N}}
\def\mm{{\cal M}}
\def\kk{{\cal K}}
\def\ll{{\cal L}}
\def\vv{{\cal V}}
\def\ff{{\cal F}}
\def\hh{{\cal H}}
\def\tt{{\cal T}}
\def\uu{{\cal U}}
\def\cc{{\cal C}}
\def\pp{{\cal P}}
\def\aa{{\cal A}}
\def\oo{{\cal O}}
\def\II{{\cal I}}
\def\dd{{\cal D}}
\def\ll{{\cal L}}
\def\ff{{\cal F}}
\def\G{{\cal G}}

\def\hs{\hat{s}}
\def\hz{\hat{z}}
\def\hL{\hat{L}}
\def\hl{\hat{l}}
\def\hl{\hat{l}}
\def\hc{\hat{\cc}}
\def\hbb{\widehat{\cal B}}
\def\hu{\hat{u}}
\def\hX{\hat{X}}
\def\hx{\hat{x}}
\def\hu{\hat{u}}
\def\hv{\hat{v}}
\def\hQ{\hat{Q}}
\def\hC{\widehat{C}}
\def\hF{\hat{F}}
\def\hf{\hat{f}}
\def\hii{\hat{\ii}}
\def\hr{\hat{r}}
\def\hq{\hat{q}}
\def\hy{\hat{y}}
\def\hZ{\widehat{Z}}
\def\hz{\hat{z}}
\def\hE{\widehat{E}}
\def\hR{\widehat{R}}
\def\hell{\hat{\ell}}
\def\hs{\hat{s}}
\def\hW{\widehat{W}}
\def\hS{\widehat{S}}
\def\hV{\widehat{V}}
\def\hB{\widehat{B}}
\def\hhh{\widehat{\cal H}}
\def\hK{\widehat{K}}
\def\hU{\widehat{U}}
\def\hhh{\widehat{\hh}}
\def\hdd{\widehat{\dd}}
\def\hZ{\widehat{Z}}

\def\hal{\hat{\alpha}}
\def\hbe{\hat{\beta}}
\def\hg{\hat{\gamma}}
\def\hrho{\hat{\rho}}
\def\hd{\hat{\delta}}
\def\hphi{\hat{\phi}}
\def\hmu{\hat{\mu}}
\def\hnu{\hat{\nu}}
\def\hsi{\hat{\sigma}}
\def\htau{\hat{\tau}}
\def\hpi{\hat{\pi}}
\def\hep{\hat{\epsilon}}
\def\hxi{\hat{\xi}}
\def\hLa{\widehat{\Lambda}^u}
\def\hPhi{\widehat{\Phi}}
\def\hPsi{\widehat{\Psi}}
\def\hPhii{\widehat{\Phi}^{(i)}}

\def\tc{\tilde{C}}
\def\tg{\tilde{\gamma}}  
\def\tV{\widetilde{V}}
\def\tC{\widetilde{\cc}}
\def\tr{\tilde{R}}
\def\tb{\tilde{b}}
\def\tt{\tilde{t}}
\def\tx{\tilde{x}}
\def\tp{\tilde{p}}
\def\tz{\tilde{Z}}
\def\tZ{\tilde{Z}}
\def\tF{\tilde{F}}
\def\tf{\tilde{f}}
\def\tp{\tilde{p}}
\def\te{\tilde{e}}
\def\tv{\tilde{v}}
\def\tu{\tilde{u}}
\def\tw{\tilde{w}}
\def\ts{\tilde{\sigma}}
\def\tr{\tilde{r}}
\def\tU{\tilde{U}}
\def\tS{\tilde{S}}
\def\tP{\widetilde{\Pi}}
\def\ttau{\tilde{\tau}}
\def\tLip{\widetilde{\Lip}}
\def\tz{\tilde{z}}
\def\tS{\tilde{S}}
\def\tts{\tilde{\sigma}}
\def\tVl{\widetilde{V}^{(\ell)}}
\def\tVj{\widetilde{V}^{(j)}}
\def\tVo{\widetilde{V}^{(1)}}
\def\tVj{\widetilde{V}^{(j)}}
\def\tPsi{\tilde{\Psi}}
 \def\tp{\tilde{p}}
 \def\tVjo{\widetilde{V}^{(j_0)}}
\def\tvj{\tilde{v}^{(j)}}
\def\tVjj{\widetilde{V}^{(j+1)}}
\def\tvl{\tilde{v}^{(\ell)}}
\def\tVt{\widetilde{V}^{(2)}}
\def\tR{\tilde{R}}
\def\tQ{\tilde{Q}}
\def\oL{\tilde{\Lambda}}
\def\tq{\tilde{q}}
\def\tx{\tilde{x}}
\def\ty{\tilde{y}}
\def\tz{\tilde{z}}
\def\txo{\tilde{x}^{(0)}}
\def\tso{\tilde{\sigma}_0}
\def\tmt{\tilde{\Lambda}}
\def\tg{\tilde{g}}
\def\tsi{\tilde{\sigma}}
\def\ttt{\tilde{t}}
\def\tC{\tilde{C}}
\def\tc{\tilde{c}}
\def\tell{\tilde{\ell}}
\def\trho{\tilde{\rho}}
\def\ts{\tilde{s}}
\def\tB{\widetilde{B}}
\def\thh{\widetilde{\cal H}}
\def\tV{\widetilde{V}}
\def\trr{\tilde{r}}
\def\te{\tilde{e}}
\def\tv{\tilde{v}}
\def\tu{\tilde{u}}
\def\tw{\tilde{w}}
\def\trho{\tilde{\rho}}
\def\tell{\tilde{\ell}}
\def\tz{\tilde{Z}}
\def\tF{\tilde{F}}
\def\tf{\tilde{f}}
\def\tp{\tilde{p}}
\def\ttau{\tilde{\tau}}
\def\tz{\tilde{z}}
\def\tg{\tilde{\gamma}}  
\def\tV{\widetilde{V}}
\def\tC{\widetilde{\cc}}
\def\tLa{\widetilde{\Lambda}^u}
\def\tR{\widetilde{R}}
\def\tr{\tilde{r}}
\def\tc{\widetilde{C}}
\def\tD{\widetilde{D}}
\def\tt{\tilde{t}}
\def\tx{\tilde{x}}
\def\tp{\tilde{p}}
\def\tS{\tilde{S}}
\def\tts{\tilde{\sigma}}
\def\tZ{\widetilde{Z}}
\def\tdelta{\tilde{\delta}}
\def\th{\tilde{h}}
\def\tB{\widetilde{B}}
\def\thh{\widetilde{\hh}}
\def\tep{\tilde{\ep}}
\def\tE{\widetilde{E}}
\def\tu{\tilde{u}}
\def\txi{\tilde{\xi}}
\def\teta{\tilde{\eta}}

\def\sr{{\sc r}}
\def\mt{{\Lambda}}
\def\do{\partial \Omega}
\def\dk{\partial K}
\def\dl{\partial L}
\def\wloc{W_{\epsilon}}
\def\piU{\pi^{(U)}}
\def\Rio{\R_{i_0}}
\def\Ri{\R_{i}}
\def\Rii{\R^{(i)}}
\def\Riii{\R^{(i-1)}}
\def\hRii{\widehat{\R}_i}
\def\hRiio{\widehat{\R}_{(i_0)}}
\def\Eii{E^{(i)}}
\def\Eio{E^{(i_0)}}
\def\Rj{\R_{j}}
\def\Vio{{\cal V}^{i_0}}
\def\Vi{{\cal V}^{i}}
\def\Wio{W^{i_0}}
\def\Wioo{W^{i_0-1}}
\def\hi{h^{(i)}}
\def\Psii{\Psi^{(i)}}
\def\pii{\pi^{(i)}}
\def\piii{\pi^{(i-1)}}
\def\gxyii{g_{x,y}^{i-1}}
\def\span{\mbox{\rm span}}
\def\Jac{\mbox{\rm Jac}}
\def\Vol{\mbox{\rm Vol}}
\def\limp{\lim_{p\to\infty}}
\def\hh{{\mathcal H}}

\def\yijl{Y_{i,j}^{(\ell)}}
\def\xijl{X_{i,j}^{(\ell)}}
\def\hyijl{\widehat{Y}_{i,j}^{(\ell)}}
\def\hxijl{\widehat{X}_{i,j}^{(\ell)}}
\def\eijl{\omega_{i,j}^{(\ell)}}
\def\J{\sf J}
\def\Gl{\Gamma_\ell}

\def\hLao{\widehat{\Lambda}^{u,1}}
\def\tLao{\widetilde{\Lambda}^{u,1}}
\def\Lao{\Lambda^{u,1}}
\def\cLao{\check{\Lambda}^{u,1}}
\def\cB{\check{B}}
\def\tpi{\tilde{\pi}}

\noindent
{\large\bf Spectra of  Ruelle transfer operators for Axiom A flows}

\bs

\noindent
{\bf Luchezar Stoyanov}

\bs

\footnotesize

\noindent
{\bf Abstract.} For Axiom A flows on basic sets satisfying certain additional conditions  
we prove  strong  spectral estimates for Ruelle transfer operators similar to these of 
Dolgopyat \cite{kn:D2} for geodesic flows on compact surfaces (for general potentials)
and transitive Anosov flows on compact manifolds  with $C^1$ jointly non-integrable 
horocycle foliations (for the Sinai-Bowen-Ruelle potential).
Here we deal with general potentials and on spaces of arbitrary dimension, although under
some geometric and regularity conditions.
As is now well known, such results have deep implications in some related 
areas, e.g. in studying analytic properties of  Ruelle zeta functions and partial differential 
operators,  closed orbit counting functions, decay of correlations  for H\"older continuous  potentials. 

\normalsize

\section{Introduction}

\renewcommand{\theequation}{\arabic{section}.\arabic{equation}}

\subsection{Introduction and main results}

Let $\phi_t : M \longrightarrow M$ be a $C^2$ Axiom A flow on a
$C^2$ complete (not necessarily compact) Riemann manifold $M$ and let $\mt$ be a basic set for $\phi_t$.
Let $\|\cdot \|$ be the {\it norm} on $T_xM$ determined by the Riemann metric on $M$ and let
$E^u(x)$ and $E^s(x)$ ($x\in \mt$)  be the tangent spaces to the strong unstable and stable manifolds 
$W^u_\ep(x)$  and $W^s_\ep(x)$, respectively (see section 2). For any $x \in \mt$, $T > 0$ and 
$\delta\in (0,\ep]$ set
$$B^u_T (x,\delta) = \{ y\in W^u_{\ep}(x) : d(\phi_t(x), \phi_t(y)) \leq \delta \: \: , 
\:\:  0 \leq t \leq T \}\;.$$

We will say that $\phi_t$ has a {\it regular distortion along unstable manifolds} over
the basic set $\mt$  if there exists a constant $\ep_0 > 0$ with the following properties:

\ms 

(a) For any  $0 < \delta \leq   \ep \leq \ep_0$ there exists a constant $R =  R (\delta , \ep) > 0$ such that 
\be
\diam( \mt \cap B^u_T(z ,\ep))   \leq R \, \diam( \mt \cap B^u_T (z , \delta))
\ee
for any $z \in \mt$ and any $T > 0$.

\ms

(b) For any $\ep \in (0,\ep_0]$ and any $\rho \in (0,1)$ there exists $\delta  \in (0,\ep]$
such that for  any $z\in \mt$ and any $T > 0$ we have
$\diam ( \mt \cap B^u_T(z ,\delta))   \leq \rho \; \diam( \mt \cap B^u_T (z , \ep))\;.$

\bs

Part (a) of the above condition resembles the Second Volume Lemma of Bowen and Ruelle \cite{kn:BR} 
about balls in Bowen's metric; this time however we deal with diameters instead of volumes.
In a separate paper \cite{kn:St3} we describe a rather general class of flows on basic sets 
satisfying this condition -- see section 7 below for a brief account of these.
There are reasons to believe that this condition may actually hold for all $C^2$
flows on basic sets -- see the comments at the end of this subsection.

In this paper we deal with flows $\phi_t$ over basic sets $\mt$ having a regular distortion 
along unstable manifolds. Apart from that, in the main result below we impose an additional condition 
on $\phi_t$ and $\mt$, called the {\it local non-integrability condition} (LNIC); it is stated in 
section  2 below. 
It should be mentioned that this condition is rather weak and is
always  satisfied e.g. for contact flows that are either Anosov 
(i.e. $\mt = M$), or have one-dimensional (un)stable manifolds (see section 7 below). 
One would  expect that (LNIC) is satisfied in most interesting cases. For example, it was shown in 
\cite{kn:St2} that open billiard flows (in any dimension) with $C^1$ (un)stable laminations over the 
non-wandering set $\mt$ always satisfy (LNIC).

Let $\rr = \{R_i\}_{i=1}^k$ be a Markov family for $\phi_t$ over $\mt$ consisting of 
rectangles $R_i = [U_i ,S_i ]$, where $U_i$ (resp. $S_i$) are (admissible) subsets of 
$W^u_{\ep}(z_i) \cap \mt$
(resp. $W^s_{\ep}(z_i) \cap \mt$) for some $\ep > 0$ and $z_i\in \mt$ (cf. section 2 for details). 
Assuming that the local stable and unstable laminations over $\mt$ are Lipschitz, the 
{\it first return time function} 
$\tau : R = \cup_{i=1}^k R_i \longrightarrow [0,\infty)$ 
and the standard Poincar\'e map $\pp: R \longrightarrow R$ are essentially Lipschitz.
Setting $U = \cup_{i=1}^k U_i$, the {\it shift map} $\sigma : U \longrightarrow U$
is defined by $\sigma = \piU\circ \pp$, where $\piU : R \longrightarrow U$
is the projection along the leaves of local stable manifolds. Let $\hU$ be the set of all
$u \in U$ whose orbits do not have common points with the boundary of $R$ (see section 2).
Given a Lipschitz real-valued function $f$  on $\hU$, set $g = g_f = f - P\tau$, where 
$P = P_f\in \R$ is the unique 
number such that the topological pressure $\Pr_\sigma(g)$ of $g$ with respect to $\sigma$ is 
zero (cf. e.g. \cite{kn:PP}). For $a, b\in \R$, one defines the {\it Ruelle transfer operator}
$L_{g-(a+\i b)\tau} : \clip (\hU) \longrightarrow \clip (\hU)$ in the usual way (cf. section 2). Here
$\clip (\hU)$ is the space of Lipschitz functions $g: \hU \longrightarrow \C$. By 
$\Lip(g)$ we denote the Lipschitz constant of $g$ and  by $\| g\|_0$ the {\it standard $\sup$ norm}  
of $g$ on $\hU$.

We will say  that the {\it Ruelle transfer operators related to the function $f$ on $\hU$ are 
eventually contracting} 
if for every $\epsilon > 0$ there exist constants $0 < \rho < 1$, $a_0 > 0$ and  $C > 0$ such 
that if $a,b\in \R$  
satisfy $|a| \leq a_0$ and $|b| \geq 1/a_0$, then for every integer $m > 0$ and every  
$h\in \clip (\hU)$ we have
$$\|L_{f -(P_f+a+ \i b)\tau}^m h \|_{\lip,b} \leq C \;\rho^m \;|b|^{\ep}\; \| h\|_{\lip,b}\; ,$$
where the norm $\|.\|_{\lip,b}$ on $\clip (\hU)$ is defined by 
$\| h\|_{\lip,b} = \|h\|_0 + \frac{\lip(h)}{|b|}$.
This implies in particular that the spectral radius  of $L_{f-(P_f+ a+\i b)\tau}$ on $\clip(\hU)$ 
does not exceed  $\rho$.

Our main result in this paper is the following.

\bs

\noindent
{\bf Theorem 1.1.} {\it  Let $\phi_t : M \longrightarrow M$ be  a $C^2$ Axiom A flow on
a $C^2$ complete  Riemann manifold satisfying the condition (LNIC) and having a
regular distortion along unstable manifolds over a  basic set $\mt$. Assume in addition  
that the local holonomy maps
along stable laminations through $\mt$ are uniformly Lipschitz. Then for any  Lipschitz real-valued function 
$f$  on $\hU$ the Ruelle transfer operators related to  $f$  are eventually contracting}.

\bs

We refer the reader to section 2 below for the definition of holonomy maps. 
In general these  are only H\"older continuous. It is known that uniform
Lipschitzness of the local stable holonomy maps can be derived from certain bunching 
condition concerning the rates of expansion/contraction of the flow along local unstable/stable 
manifolds over $\mt$ (see \cite{kn:Ha1}, \cite{kn:Ha2}, \cite{kn:PSW}). 

It should be mentioned that some kind of a non-integrability assumption is necessary for results like 
Theorem 1.1 (and  Corollary 1.2). Indeed, Pollicott \cite{kn:Po} and Ruelle \cite{kn:R2} have 
constructed examples of mixing  Axiom A flows with jointly integrable stable and
unstable laminations for which the statements of Corollaries 1.4 and 1.5 below (and therefore that of 
Theorem 1.1 as well) do not hold. As one can see in section 2 below, (LNIC) is a rather weak 
non-integrability condition. Moreover, in the contact case it follows from another
condition (see (ND) in section 6 below) which looks rather natural and is always satisfied for 
Anosov flows.

It follows from Theorem 7.1  below that a flow with one-dimensional unstable laminations over a 
basic set $\mt$ always has a regular distortion along unstable manifolds. The local stable 
holonomy maps are always Lipschitz (in fact $C^1$) if the stable laminations over $\mt$ are 
one-dimensional (see e.g. Theorem 1 and fact (2) on p. 647 in \cite{kn:Ha1}).
Moreover, the flow always satisfies (LNIC) if it is Anosov with jointly non-integrable stable and unstable
foliations, or it is contact (see Proposition 6.1 for the latter). 
Thus, we get the  following consequence of Theorem 1.1.

\bs

\noindent
{\bf Corollary 1.2.} {\it Let $\phi_t : M \longrightarrow M$ be a $C^2$ flow on a $C^2$ Riemann manifold 
and let $\mt$ be a  basic set for $\phi_t$ such that the stable and unstable laminations over $\mt$ are 
one-dimensional. Assume in addition one of the following: (i) the flow is contact, or (ii) the flow is 
Anosov  and the stable and unstable laminations over $\mt$ are jointly non-integrable.  
Then the Ruelle transfer  operators related to any Lipschitz real-valued function $f$  on $\hU$ 
are eventually contracting.}


\bs

The above  was first proved by Dolgopyat (\cite{kn:D1}, \cite{kn:D2}) in the case of geodesic flows on 
compact surfaces.
The second main result in \cite{kn:D2} concerns transitive Anosov flows  on
compact Riemann manifolds with $C^1$ jointly non-integrable local stable and unstable foliations.
For such flows Dolgopyat proved that the conclusion of Theorem 1.1 holds for the 
Sinai-Bowen-Ruelle potential  $f = \log \det (d\phi_\tau)_{|E^u}$.
Theorem 1.1 appears to be first result of this kind that works for any potential and in any 
dimension\footnote{Albeit under additional conditions but 
this appears to be inevitable.}.
  
For contact Anosov flows in any dimension (LNIC) is always satisfied (see Proposition 6.1 below), 
so the following is also an immediate consequence  of Theorem 1.1.

\bs

\noindent
{\bf Corollary 1.3.} {\it  Let $\phi_t : M \longrightarrow M$ be a $C^2$ contact Anosov flow on
a compact Riemann manifold  having a regular distortion along unstable manifolds over  
$M$ and such that the local holonomy maps along stable laminations through $\mt$ are uniformly 
Lipschitz.  Then the Ruelle transfer operators related to any Lipschitz real-valued function  
$f$  on $\hU$ are eventually contracting.}

\bs

Using a smoothing procedure as in \cite{kn:D2} (see also Corollary 3.3 in \cite{kn:St1}), an estimate
similar to that in Theorem 1.1  holds for the Ruelle operator acting on the space 
$\ff_\gamma(U)$ of H\"older continuous functions with respect to an appropriate norm $\|\cdot \|_{\gamma,b,U}$.

Using Theorem 1.1  and an argument of Pollicott and Sharp \cite{kn:PoS1} 
one derives valuable information  about the {\it Ruelle zeta function}
$\zeta(s) = \prod_{\gamma} (1- e^{-s\ell(\gamma)})^{-1}\;,$
where $\gamma$ runs over the set of primitive  closed orbits of $\phi_t: \mt \longrightarrow \mt$
and $\ell(\gamma)$ is the least period of $\gamma$.  In what follows $h_T$ denotes the
 {\it topological entropy} of $\phi_t$ on $\mt$.\\

\noindent
{\bf Corollary 1.4.} {\it Under the assumptions in Theorem 1.1, Corollary 1.2 or Corollary 1.3, 
the  zeta function $\zeta(s)$ of the flow $\phi_t: \mt \longrightarrow \mt$ has an analytic  
and non-vanishing continuation in a half-plane $\Re(s) > c_0$ for some $c_0 < h_T$ except for a 
simple pole at $s = h_T$.  Moreover, there exists $c \in (0, h_T)$ such that
$$\pi(\lambda) = \# \{ \gamma : \ell(\gamma) \leq \lambda\} = \li(e^{h_T \lambda}) + O(e^{c\lambda})$$
as $\lambda\to \infty$, where $\di \li(x) = \int_2^x \frac{du}{\log u} \sim \frac{x}{\log x}$ as 
$x \to \infty$. }

\bigskip

As another consequence of Theorem 1.1  and the procedure described in \cite{kn:D2} one gets exponential 
decay of correlations for the  flow $\phi_t : \mt \longrightarrow \mt$. 

Given $\alpha > 0$ denote by $\ff_\alpha(\mt)$ the set of {\it H\"older  continuous functions} with 
H\"older exponent $\alpha$ and by $\| h\|_\alpha$ the {\it H\"older constant} of $h\in \ff_\alpha(\mt)$. \\

\noindent
{\bf Corollary 1.5.} {\it   Under the assumptions in Theorem 1.1, Corollary 1.2  or Corollary 1.3,
let $F$ be a H\"older continuous function on $\mt$  and let $\nu_F$ be the Gibbs measure
determined by $F$ on $\mt$. Assume in addition that the manifold $M$ and the flow $\phi_t$ are $C^5$.
For every $\alpha > 0$ there exist constants $C = C(\alpha) > 0$ and 
$c = c(\alpha) > 0$ such that 
$$\left| \int_{\mt} A(x) B(\phi_t(x))\; d\nu_F(x) - 
\left( \int_{\mt} A(x)\; d\nu_F(x)\right)\left(\int_{\mt} B(x) \; d\nu_F(x)\right)\right|
\leq C e^{-ct} \|A\|_\alpha \; \|B\|_\alpha \;$$
for any two functions $A, B\in \ff_\alpha(\mt)$.}

\bs

There has been a considerable activity in recent times to establish exponential and other types of decay of
correlations for various  kinds of systems with some highly rated results of Chernov \cite{kn:Ch1},
Dolgopyat \cite{kn:D1}, \cite{kn:D2}, Liverani \cite{kn:L1}, \cite{kn:L2}, Young \cite{kn:Y1}, \cite{kn:Y2}.
See also \cite{kn:BSC}, \cite{kn:BaT}, \cite{kn:T}, \cite{kn:Ch2}, \cite{kn:Ch3}, \cite{kn:GL}, 
\cite{kn:D3},  \cite{kn:FMT}, \cite{kn:Mel},  \cite{kn:ChY} and the references there. 
In \cite{kn:L2} Liverani proves exponential decay of correlations for contact Anosov flows, and this 
appears to be the most general result of this kind  so far. Recently Tsujii \cite{kn:T}
obtained finer results for the same kind of flows. It should be stressed that in \cite{kn:L2} 
and \cite{kn:T} (and various other works; see the references there) spectral properties of  a 
different kind of transfer operators are studied, namely the operators 
$\ll_tg = g\circ \phi_{-t}$ ($t\in \R$) acting on functions $g$ on a compact manifold $M$, 
$\phi_t$ being a contact Anosov flow on $M$. 

We should stress though that the {\bf main aim of this article is not to establish results on decay
of correlations but rather to get strong spectral estimates for Ruelle transfer operators}.
These operators appear to be more difficult to deal with, 
since they are intimately related to geometric properties of the flow (and the basic set $\mt$ 
when Axiom A flows are considered). On the other hand, spectral results of the kind obtained in 
Theorem 1.1 appear to be much finer and to have a wider and deeper range of applications.

In \cite{kn:St1} a modification of the method from \cite{kn:D2} was used to prove an analogue
of Corollary 1.2 above for open billiard flows in the plane. Using similar tools, Naud \cite{kn:N}
proved a similar result for geodesic flows on convex co-compact hyperbolic surfaces. The results in
both \cite{kn:St1} and \cite{kn:N} are special cases of Corollary 1.2.
Baladi and Vall\'ee (\cite{kn:BaV})  obtained Dolgopyat type estimates for 
transfer operators in the case of a suspension of an interval map.

It has been well known since Dolgopyat's paper \cite{kn:D2} that
strong spectral estimates for Ruelle transfer operators as the ones described  in Theorem 1.1 
lead to deep results concerning zeta functions and related topics which are difficult to obtain by 
other means. For example, such estimates were fundamental in \cite{kn:PoS1}, where the statement of
Corollary 1.4 was proved for geodesic flows on compact surfaces of negative curvature. For the same kind 
of flows, fine and very interesting asymptotic estimates for pairs of closed geodesics were established 
in \cite{kn:PoS3}, again by using the strong spectral estimates in \cite{kn:D2}. For Anosov flows 
with $C^1$  jointly non-integrable horocycle foliations full asymptotic expansions for counting 
functions similar to $\pi(\lambda)$ however with some homological constraints were obtained in 
\cite{kn:An} and \cite{kn:PoS2}. In \cite{kn:PeS2} Theorem 1.1 above was used to obtain results similar 
to these in \cite{kn:PoS3} about correlations 
for pairs of closed billiard trajectories  for billiard flows in $\R^n \setminus K$,
where $K$ is a finite disjoint union of strictly convex compact bodies with smooth boundaries 
satisfying the so  called `no eclipse condition'  (and some additional conditions as well). 
For the same kind of models and using  Theorem 1.1 again, a rather non-trivial result was 
established in \cite{kn:PeS1}  about analytic  continuation of the cut-off resolvent of the 
Dirichlet Laplacian in $\R^n\setminus K$, which appears to be the first of its kind in 
the field of quantum chaotic scattering. It is not clear at all how such a result could be 
proved without using the strong spectral estimates of the kind considered here. Finally, in a very
recent preprint \cite{kn:PeS3}, using  Theorem 1.1 and under the assumptions in this theorem,
a fine asymptotic was obtained  for the number of closed trajectories in $\mt$
with primitive periods lying in exponentially shrinking intervals 
$(x - e^{-\delta x}, x + e^{-\delta x})$, $\delta > 0$, $x \to + \infty.$

The main part of this article consists of sections 3-5 where Theorem 1.1 is proved.
Section 2 contains some basic definitions and  facts. The regular distortion along  unstable manifolds 
has some natural consequences about diameters of cylinders defined by means of Markov families --
these are derived in section 3. The essential part of the proof of Theorem 1.1 is in sections 4-5. 

In section 6 we introduce a non-integrability condition for contact flows and show that, under some
regularity assumption, it implies (LNIC).
In section 7 we describe classes of flows over basic sets satisfying the conditions in Theorem 1.1
using the results in \cite{kn:St3}.
It should be stressed that the central  part of the arguments in \cite{kn:St3} is to prove a local 
version of regular distortion along unstable manifolds, where  e.g. (1.1) is satisfied at a single point 
$z\in \mt$ (with a constant $R = R(z,\delta, \ep)$ depending on $z$ as well). 
It is not difficult to see (using a variation of the arguments 
in \cite{kn:St3}) that a similar local result can be proved at Lyapunov regular points $z\in \mt$ for 
any $C^2$ flow on any basic set\footnote{Now the role of the 'bottom of the unstable spectrum' is 
played by the exponential of the least positive Lyapunov exponent.}. 
It is not yet clear how one can get a uniform global result over $\mt$ from such local  results.

The Appendix  contains the proof of  a technical lemma from section 5.

\subsection{Comments on the proof of the main result}

In the proof of Theorem 1.1 we use the general framework of Dolgopyat's method from \cite{kn:D1}, \cite{kn:D2} 
and its modification in \cite{kn:St1}, however a significant  new development is necessary. 
The main difficulty comes from the fact that very little is known about the (geometric) structure 
of the set $\mt$. For general potentials $f$ the original method in \cite{kn:D2} was only applied 
to geodesic flows on surfaces. It requires more sophisticated arguments to deal with general 
potentials for flows on higher-dimensional spaces, especially when the flow is considered over a 
basic set (possibly with a complicated fractal structure)
rather than on a nice smooth compact manifold.

In this subsection we make some general remarks about the main points in the proof of Theorem 1.1.
Given $f \in \clip(U)$, one wants to show that large powers of the Ruelle transfer operator
$L_{r} = L_{f-(P_f+a+\i \, b)\tau}$ are contracting for small $a\in \R$ and large $b\in \R$. 
For any integer $N > 0$ we have
\be
L^N_{r} h(x) = \sum_{\sigma^N(y) = x} e^{r_N(y)}\, h(y)\;,
\ee
where $r_N(y) = r(y) + r(\sigma(y)) + \ldots + r(\sigma^{N-1}(y))$. 
One of the main steps in \cite{kn:D2} is to define appropriately $C^1$ inverses 
$v_1, v_2 : U_0\longrightarrow U$  of $\sigma^N$, i.e. $\sigma^N(v_i(x)) = x$ for $x$ in some `small' 
open subset $U_0$ of $U$. For these one ultimetely shows that for large $N$ there exists $\lambda\in (0,1)$ 
such that
\begin{eqnarray}
\qquad \quad |e^{r_N(v_1(x))}\, h(v_1(x)) + e^{r_N(v_2(x))}\, h(v_2(x)) |
 \leq \lambda\, [ |e^{r_N(v_1(x))}|\, h(v_1(x)) +  | e^{r_N(v_2(x))} |\, h(v_2(x)) ]
\end{eqnarray}
for small  $|a|$ and large $|b|$ (of magnitude depending on $N$) and $h\in C^1(U)$ with $h > 0$, $\|h\| \leq \Con$,
$\|dh\|\leq \Con\, |b|$. This leads to a similar estimate for the whole sum in (1.2).  To make this cancellation 
mechanism\footnote{To use the words of Liverani \cite{kn:L2}
who is using a similar idea.} work  for general complex-valued functions $h$,  one uses estimates 
involving $\nn_J |h|$, where $\{ \nn_J \}$ is a finite family of specially defined  operators acting on positive functions 
in $C^1(U)$ with bounded logarithmic derivatives. One of the main features of Dolgopyat's operators
$\nn_J$ is that they are $L^2$-contractions with respect to the invariant Gibbs measure $\nu$ determined by
$f- P\, \tau$ on $U$, and this is crucial for the proof of  the main result.

The whole procedure is rather more complicated than the above, and we refer the reader to 
section 5 below for more details. It is worth mentioning though that it is the construction of the inverses 
$v_1$ and $v_2$ of $\sigma^N$ and the proof of their main properties where the joint non-integrability of the 
stable and unstable families is used. In \cite{kn:D2} this  results in finding $C^1$ vector fields 
$e_1(z), \ldots, e_n(z)$ defined in a small neighbourhood $U_0$ of a point $z_0 \in U$ so that
for large $N$,
\be
|\partial _{e_1}(\tau_N(v_1(u)) - \tau_N(v_2(u)))| \geq \ep \quad, \quad u\in U_0\;,
\ee
\be
|\partial _{e_i}(\tau_N(v_1(u)) - \tau_N(v_2(u)))| << \ep \:\:\:\:, \: i > 1\;,\: u \in U_0\;.
\ee
Notice that $e^{\i\, b\, \tau_N(y)}$ is the `tricky' part of the exponential term $e^{r_N(y)}$ in (1.2).
With (1.4) and (1.5) one has that if $u,u'\in U_0$ are $C \ep$-close however there is a $c \ep$-gap 
between their first coordinates (with respect to the vector field $e_1$) for some constants $C>  c > 0$, then 
\be
\left|\left[\tau_N(v_1(u)) - \tau_N(v_2(u))\right] - \left[\tau_N(v_1(u')) - \tau_N(v_2(u'))\right]\right|
\geq \con\, \ep\;.
\ee
This is what lies beneath the proof of (1.3).

In \cite{kn:St1} a modification of the above was used to deal with open billiard flows in the plane. In this case the unstable 
manifolds are  one-dimensional, so one has just one vector field $e_1(u)$, however $U$ is a Cantor set, so the construction of 
$v_1(u)$ and $v_2(u)$ is non-trivial.  One of the main difficulties in \cite{kn:St1} was to show (using specific features of the model) 
that there exist constants $c_2 > c_1 > 0$ such that for any $\ep > 0$, $U$ can be partitioned into intervals of lengths between 
$c_1\ep$ and $c_2\ep$ such that successive intervals with common points with $\mt$ have a gap of a similar size between 
them\footnote{The same idea was later used in \cite{kn:N} to deal with limit sets of convex co-compact hyperbolic surfaces.}. 
Then the one-dimensionality of the unstable manifolds allows to construct  the functions $v_i(u)$ and prove (1.4). Moreover for 
any of the small intervals $\Delta_1$ in the partition of $U$ intersecting 
$\mt$ one can find another interval $\Delta_2$ intersecting $\mt$ with a $c_1\ep$-gap between $\Delta_1$ and $\Delta_2$
such that (1.6) holds for all $u \in \Delta_1$ and all $u'\in \Delta_2$. The rest of Dolgopyat's method was not so difficult to apply, 
although extra modifications were necessary to deal with the singularity of the Gibbs measures in  this case.

In the present paper the dimension $n$ of the unstable manifolds is arbitrary and the set $\mt$ 
can be `anything from a Cantor set to a manifold'.
Therefore to construct vector fields $e_i(u)$ in a neighbourhood $U_0$ of a point $z_0\in U$
with (1.4) and (1.5)  is meaningless, unless we know that there are `many points' of $\mt$ in the direction
of $e_1(u)$ (or `very close' to it) for a `large set' of $u$'s in $U_0$. To establish something 
like this however 
seems impossible without assuming anything in the spirit of the extra condition (LNIC) from 
section 2 below.  With (LNIC),  one 
constructs for large $N > 0$ a pair of  inverses $v_1(u)$ and $v_2(u)$  of $\sigma^N$ such that the 
analogue of (1.4) holds on a neighbourhood of $z$ with $e_1$ replaced by unit vectors $\eta$ in a 
cone whose axis is a parallel translate of $\eta$.  
The above gives that if $\cc_1$ and $\cc_2$ are two
cylinders  in the vicinity of $z$ of size $\leq \delta$,  for some small $\delta > 0$, that can be `separated 
by a plane whose normal  is  close to $\eta$' (one can make sense of this by using  some
parametrization of $U$ near $z$), then (1.6) holds for all $u \in \cc_1$ and all $u'\in \cc_2$. 
However, since $n > 1$ in general, just one direction $\eta$ does not give enough opportunities
to separate cylinders.  So, one needs to construct  a finite number of 
directions $\eta_1, \ldots, \eta_{\ell_0}$  tangent to  $\mt$ (at various points close to the initial point $z_0$), and for each
$j = 1, \ldots, \ell_0$, a pair of inverses $v_1^{(j)}(u)$ and $v_2^{(j)}(u)$ of $\sigma^N$ defined
on some small open subset $U_0$ of $U$ such that `sufficiently many'
pairs of cylinders in $U_0$ of size $\leq \delta$ can be separated by  planes each of them having a
'normal   close to some $\eta_j$', and this can be done for all $\delta$ in some interval 
$(0,\delta']$. This is the content of the main Lemma 4.2.

From Lemma 4.2 and the other constructions in section 4  one  gets some
analogue of (1.6) involving different pairs of inverses $v_1^{(\ell)}(u)$ and $v_2^{(\ell)}(u)$, 
and this turns out
to be enough to implement an essentially modified  analytic part of Dolgopyat's method. This is done
in section 5. Significant new development is necessary due again to the unknown structure of  
$\mt$. What saves a lot of potential
extra problems is the fact that we work with a new metric $D$ on $U$ (or rather on the subset $\hU$ of $U$) 
defined by means of cylinders. It turns out that this metric fits well with the (modified)
Dolgopyat operators $\nn_J$ and makes them work in a truly multidimensional situation. 
Moreover, with this new metric there is no
need for the Gibbs measure $\nu$ to have the so called Federer property (see \cite{kn:D2}). 
In fact, it is not 
clear at all whether $\nu$ has this property  in the cases we consider.

\def\Intu{\Int^u}
\def\Ints{\Int^s}

\section{Preliminaries}
\setcounter{equation}{0}

Throughout this paper $M$ denotes a $C^2$ complete (not necessarily compact) 
Riemann manifold,  and $\phi_t : M \longrightarrow M$ ($t\in \R$) a $C^2$ flow on $M$. A
$\phi_t$-invariant closed subset $\mt$ of $M$ is called {\it hyperbolic} if $\mt$ contains
no fixed points  and there exist  constants $C > 0$ and $0 < \lambda < 1$ such that 
there exists a $d\phi_t$-invariant decomposition 
$T_xM = E^0(x) \oplus E^u(x) \oplus E^s(x)$ of $T_xM$ ($x \in \mt$) into a direct 
sum of non-zero linear subspaces,
where $E^0(x)$ is the one-dimensional subspace determined by the direction of the flow
at $x$, $\| d\phi_t(u)\| \leq C\, \lambda^t\, \|u\|$ for all  $u\in E^s(x)$ and $t\geq 0$, and
$\| d\phi_t(u)\| \leq C\, \lambda^{-t}\, \|u\|$ for all $u\in E^u(x)$ and  $t\leq 0$.

The flow $\phi_t$ is called an {\it Axiom A flow} on $M$ if the non-wandering set 
of $\phi_t$  is a disjoint union of a finite set  consisting of  fixed
hyperbolic points and a compact hyperbolic subset containing no fixed points in which
the periodic points are dense (see e.g. \cite{kn:KH}). 
A non-empty compact $\phi_t$-invariant hyperbolic subset $\mt$ of $M$ which is not a single 
closed orbit is called a {\it basic set} for $\phi_t$ if $\phi_t$ is transitive on $\mt$ 
and $\mt$ is locally maximal, i.e. there exists an open neighbourhood $V$ of
$\mt$ in $M$ such that $\mt = \cap_{t\in \R} \phi_t(V)$. When $M$ is compact and $M$ itself
is a basic set, $\phi_t$ is called an {\it Anosov flow}.

For $x\in \Lambda$ and a sufficiently small $\epsilon > 0$ let 
$$\wloc^s(x) = \{ y\in M : d (\phi_t(x),\phi_t(y)) \leq \epsilon \: \mbox{\rm for all }
\: t \geq 0 \; , \: d (\phi_t(x),\phi_t(y)) \to_{t\to \infty} 0\: \}\; ,$$
$$\wloc^u(x) = \{ y\in M : d (\phi_t(x),\phi_t(y)) \leq \epsilon \: \mbox{\rm for all }
\: t \leq 0 \; , \: d (\phi_t(x),\phi_t(y)) \to_{t\to -\infty} 0\: \}$$
be the (strong) {\it stable} and {\it unstable manifolds} of size $\epsilon$. Then
$E^u(x) = T_x \wloc^u(x)$ and $E^s(x) = T_x \wloc^s(x)$. 
Given $\delta > 0$, set $E^u(x;\delta) = \{ u\in E^u(x) : \|u\| \leq \delta\}$;
$E^s(x;\delta)$ is defined similarly. 

It follows from the hyperbolicity of $\mt$  that if  $\epsilon_0 > 0$ is sufficiently small,
there exists $\ep_1 > 0$ such that if $x,y\in \mt$ and $d (x,y) < \ep_1$, 
then $W^s_{\ep_0}(x)$ and $\phi_{[-\ep_0,\ep_0]}(W^u_{\ep_0}(y))$ intersect at exactly 
one point $[x,y ] \in \mt$  (cf. \cite{kn:KH}). That is, there exists a unique 
$t\in [-\ep_0, \ep_0]$ such that
$\phi_t([x,y]) \in W^u_{\ep_0}(y)$. Setting $\Delta(x,y) = t$, 
defines the so called {\it temporal distance
function} (\cite{kn:KB},\cite{kn:Ch1}, \cite{kn:D2}) which will be used significantly
throughout this paper. 
For $x, y\in \mt$ with $d (x,y) < \ep_1$, define
$\pi_y(x) = [x,y] = W^s_{\ep}(x) \cap \phi_{[-\ep_0,\ep_0]} (W^u_{\ep_0}(y))\;.$
Thus, for a fixed $y \in \mt$, $\pi_y : W \longrightarrow \phi_{[-\ep_0,\ep_0]} (W^u_{\ep_0}(y))$ is the
{\it projection} along local stable manifolds defined on a small open neighbourhood $W$ of $y$ in $\mt$.
Choosing $\ep_1 \in (0,\ep_0)$ sufficiently small, 
the restriction 
$\pi_y: \phi_{[-\ep_1,\ep_1]} (W^u_{\ep_1}(x)) \longrightarrow \phi_{[-\ep_0,\ep_0]} (W^u_{\ep_0}(y))$
is called a {\it local stable holonomy map\footnote{In a similar way one can define
holonomy maps between any two sufficiently close local transversals to stable laminations; see e.g.
\cite{kn:PSW}.}.} Combining such a map with a shift along the flow we get another local stable holonomy  map
$\hh_x^y : W^u_{\ep_1}(x) \cap \mt  \longrightarrow W^u_{\ep_0}(y) \cap \mt$.
In a similar way one defines local holonomy maps along unstable laminations.

Given $z \in \mt$, let $\exp^u_z : E^u(z;\ep_0) \longrightarrow W^u_{\ep_0}(z)$  and
$\exp^s_z : E^s(z;\ep_0) \longrightarrow W^s_{\ep_0}(z)$ be the corresponding
{\it exponential maps}.
A  vector $\eta\in E^u(z)\setminus \{ 0\}$ will be called  {\it tangent to $\mt$} at
$z$ if there exist infinite sequences $\{ v^{(m)}\} \subset  E^u(z)$ and $\{ t_m\}\subset \R\setminus \{0\}$
such that $\exp^u_z(t_m\, v^{(m)}) \in \mt \cap W^u_{\ep}(z)$ for all $m$, $v^{(m)} \to \eta$ and 
$t_m \to 0$ as $m \to \infty$. 
It is easy to see that a vector $\eta\in E^u(z)\setminus \{ 0\}$ is  tangent to $\mt$ at
$z$ if there exists a $C^1$ curve $z(t)$ ($0\leq t \leq a$) in $W^u_{\ep}(z)$ for some $a > 0$ 
with $z(0) = z$ and $\dot{z}(0) = \eta$ such that $z(t) \in \mt$ for arbitrarily small $t >0$.

The following is the  {\it local non-integrability condition}\footnote{In a previous version of this
paper,  we required (2.1) with $\ty_2 = \tz$, i.e. 
$|\Delta( \exp^u_{z}(v), \pi_{\ty_1}(z))| \geq \delta\,  \|v\|$. 
The present version of (LNIC) appears to be significantly weaker. See the Remark after
the condition (ND) in section 6.} for $\phi_t$ and $\mt$ mentioned in section 1.

\medskip

\noindent
{\sc (LNIC):}  {\it There exist $z_0\in \mt$,  $\ep_0 > 0$ and $\theta_0 > 0$ such that
for any  $\ep \in (0,\ep_0]$, any $\hz\in \mt \cap W^u_{\ep}(z_0)$  and any tangent vector 
$\eta \in E^u(\hz)$ to $\mt$ at $\hz$ with 
$\|\eta\| = 1$ there exist  $\tz \in \mt \cap W^u_{\ep}(\hz)$, 
$\ty_1, \ty_2 \in \mt \cap W^s_{\ep}(\tz)$ with $\ty_1 \neq \ty_2$,
$\delta = \delta(\tz,\ty_1, \ty_2) > 0$ and $\ep'= \ep'(\tz,\ty_1,\ty_2)  \in (0,\ep]$ such that
\be
|\Delta( \exp^u_{z}(v), \pi_{\ty_1}(z)) -  \Delta( \exp^u_{z}(v), \pi_{\ty_2}(z))| \geq \delta\,  \|v\| 
\ee
for all $z\in W^u_{\ep'}(\tz)\cap\mt$  and  $v\in E^u(z; \ep')$ with  $\exp^u_z(v) \in \mt$ and
$\la \frac{v}{\|v\|} , \eta_z\ra \geq \theta_0$,   where $\eta_z$ is the parallel 
translate of $\eta$ along the geodesic in $W^u_{\ep_0}(z_0)$ from $\hz$ to $z$. }

\medskip
See Figure 2 on p. 32.
It should be mentioned that if in (LNIC) one requires $\ty_2 = \tz$, this would replace (2.1) by 
$|\Delta( \exp^u_{z}(v), \pi_{\ty}(z))| \geq \delta\,  \|v\| $ with $\ty = \ty_1$, which is still a rather 
general non-integrability condition. However in its present form
(LNIC) is a substantially weaker condition.
It is easy to see that the uniform non-integrability condition
(UNI) of Chernov \cite{kn:Ch1} and Dolgopyat  \cite{kn:D2} implies (LNIC)\footnote{In fact,
Chernov and Dolgopyat used (UNI) only for Anosov flows on 3-dimensional
manifolds. It is quite clear that when $\dim E^u(x) > 1$ ($x\in \mt$), (LNIC) is a much weaker 
condition than (UNI).}.

We will say that $A$ is an {\it admissible subset} of $W^u_{\ep}(z) \cap \mt$ ($z\in \mt$) 
if $A$ coincides with the closure
of its interior in $W^u_\ep(z) \cap \mt$. Admissible subsets of $W^s_\ep(z) \cap \mt$ are defined similarly.
Following \cite{kn:D2}, a subset $R$ of $\mt$ will be called a {\it rectangle} if it has the form
$R = [U,S] = \{ [x,y] : x\in U, y\in S\}$, where $U$ and $S$ are admissible 
subsets of $W^u_\ep(z) \cap \mt$ and $W^s_\ep(z) \cap \mt$, respectively, for some $z\in \mt$. 
In what follows we will denote by 
$\Intu(U)$ the {\it interior} of $U$ in the set $\wloc^u(z) \cap \mt$. In a similar way we 
define $\Ints(S)$, and then set $\Int(R) = [\Intu(U), \Ints(S)]$.
Given $\xi = [x,y] \in R$, set $W^u_R(\xi) = [U,y] = \{ [x',y] : x'\in U\}$ and
$W^s_R(\xi) = [x,S] = \{[x,y'] : y'\in S\} \subset W^s_{\ep_0}(x)$. The {\it interiors}
of these sets in the corresponding leaves are defined by $\Intu(W^u_R(\xi)) = [\Intu(U),y]$
and $\Ints(W^s_R(\xi)) = [x,\Ints(S)]$.

Let $\rr = \{ R_i\}_{i=1}^k$ be a family of rectangles with $R_i = [U_i  , S_i ]$,
$U_i \subset \wloc^u(z_i) \cap \mt$ and $S_i \subset \wloc^s(z_i)\cap \mt$, respectively, 
for some $z_i\in \mt$. 
Set $R =  \cup_{i=1}^k R_i$.
The family $\rr$ is called {\it complete} if  there exists $T > 0$ such that for every $x \in \mt$,
$\phi_{t}(x) \in R$ for some  $t \in (0,T]$.  Given a complete family $\rr$, the related
{\it Poincar\'e map} $\pp: R \longrightarrow R$
is defined by $\pp(x) = \phi_{\tau(x)}(x) \in R$, where
$\tau(x) > 0$ is the smallest positive time with $\phi_{\tau(x)}(x) \in R$.
The function $\tau$  is called the {\it first return time}  associated with $\rr$. 
A complete family $\rr = \{ R_i\}_{i=1}^k$ of rectangles in $\mt$ is called a 
{\it Markov family} of size $\chi > 0$ for the  flow $\phi_t$ if $\diam(R_i) < \chi$ for all $i$ and: 
(a)  for any $i\neq j$ and any $x\in \Int(R_i) \cap \pp^{-1}(\Int(R_j))$ we have   
$\pp(\Ints (W_{R_i}^s(x)) ) \subset \Ints (W_{R_j}^s(\pp(x)))$ and 
$\pp(\Intu(W_{R_i}^u(x))) \supset \Intu(W_{R_j}^u(\pp(x)))$;
(b) for any $i\neq j$ at least one of the sets $R_i \cap \phi_{[0,\chi]}(R_j)$ and
$R_j \cap \phi_{[0,\chi]}(R_i)$ is empty.

The existence of a Markov family $\rr$ of an arbitrarily small size $\chi > 0$ for $\phi_t$
follows from the construction of Bowen \cite{kn:B} (cf. also  Ratner \cite{kn:Ra}). 

From now on we will assume that $\rr = \{ R_i\}_{i=1}^k$ is a fixed Markov family for  $\phi_t$
of size $\chi < \ep_0/2 < 1$. Set  $U = \cup_{i=1}^k U_i$ and $\Intu (U) = \cup_{j=1}^k \Intu(U_j)$.
The {\it shift map} $\sigma : U   \longrightarrow U$ is given by
$\sigma  = \piU \circ \pp$, where $\piU : R \longrightarrow U$ is the {\it projection} along stable leaves. 
Notice that  $\tau$ is constant on each stable leaf $W_{R_i}^s(x) = W^s_{\ep_0}(x) \cap R_i$. 
For any integer $m \geq 1$
and any function $h : U \longrightarrow \C$ define $h_m : U \longrightarrow \C$ by
$h_m(u) = h(u) + h(\sigma(u)) + \ldots + h(\sigma^{m-1}(u))$.
 
Denote by $\widehat{U}$ the {\it core} of  $U$, i.e. the set of those $x\in U$ such that 
$\pp^m(x) \in \Int(R) = \cup_{i=1}^k \Int(R_i)$ 
for all $m \in \Z$. It is well-known (see \cite{kn:B}) that $\hU$ is a residual subset of $U$ and has full
measure with respect to any Gibbs measure on $U$.
Clearly in general $\tau$ is not continuous on $U$, however $\tau$ is {\it essentially Lipschitz} on $U$ 
in the sense that there exists a constant $L > 0$ such that if $x,y \in U_i \cap \sigma^{-1}(U_j)$ 
for some $i,j$, then $|\tau(x) - \tau(y)| \leq L\, d(x,y)$.
The same applies to $\sigma : U \longrightarrow U$.  Throughout we will mainly 
work with the restrictions of $\tau$ and $\sigma$ to $\hU$. Set $\hU_i = U_i \cap \hU$.

Let $B(\hU)$ be the {\it space of  bounded functions} $g : \hU \longrightarrow \C$ with its standard norm  
$\|g\|_0 = \sup_{x\in \hU} |g(x)|$. Given a function $g \in B(\hU)$, the  {\it Ruelle transfer operator } 
$L_g : B(\hU) \longrightarrow B(\hU)$ is defined by $\di (L_gh)(u) = \sum_{\sigma(v) = u} e^{g(v)} h(v)\;.$
If $g \in B(\hU)$ is Lipschitz on $\hU$, then  $L_g$ preserves the space $\clip(\hU)$
of {\it Lipschitz functions} $g: \hU \longrightarrow \C$.

In what follows we will assume that $f$ is a {\bf fixed  real-valued function in} $\clip(\hU)$.
Let  $P = P_f$ be the unique real number so that $\Pr_\sigma(f- P\, \tau) = 0$, where
$\Pr_\sigma(h)$ is the {\it topological pressure} of $h$ with respect to the shift map $\sigma$
(see e.g.  \cite{kn:PP}). Set $g = g_f = f - P\, \tau$.

By Ruelle-Perron-Frobenius' Theorem (see e.g. chapter 2 in \cite{kn:PP}) for any real number $a$  
with $|a|$ sufficiently small,
as an operator on $\clip(\hU)$, $L_{f-(P+a)\tau}$ has a {\it largest eigenvalue} 
$\lambda_{a}$ and there exists a
 (unique) regular probability measure $\hnu_a$ on $U$ with 
 $L_{f-(P+a)\tau}^*\hnu_a = \lambda_a\, \hnu_a$, i.e.
$\int L_{f-(P+a)\tau} H \, d\hnu_a = \lambda_a\, \int H\ d\hnu_a$ for any $H \in C(U)$.
Fix a corresponding (positive) eigenfunction $h_{a} \in \clip(\hU)$ such that $\int h_{a} \, d\hnu_a = 1$. 
Then $d\nu = h_0\, d\hnu_0$ defines a {\it $\sigma$-invariant  probability measure} $\nu$ on $U$.
Since $\Pr_\sigma (f- P\tau) = 0$, it follows from the main properties of pressure (cf. e.g. chapter 3 in \cite{kn:PP}) 
that $|\Pr_\sigma(f - (P+a)\tau)| \leq |\tau|_0 \, |a|$.  Moreover, for small $|a|$ the maximal eigenvalue
$\lambda_{a}$ and the eigenfunction $h_{a}$ are Lipschitz in  $a$, so
there exist constants $a'_0 > 0$ and $C_0 > 0$ such that $|h_{a} - h_0| \leq C_0|a|$ on $\hU$ and
$|\lambda_a - 1| \leq C_0 |a|$ for  $|a| \leq a'_0$.

For $|a|\leq a'_0$,  as in \cite{kn:D2}, consider the function
$$\fa(u) = f (u) - (P+a) \tau(u) + \ln h_{a}(u) -  \ln h_{a}(\sigma(u)) - \ln \lambda_{a}\;$$
and the operators 
$$ \lab = L_{\fa - \i\,b\tau} : \clip(\hU) \longrightarrow \clip(\hU)\:\:\: , \:\:\:
\ma = L_{\fa} : \clip(\hU) \longrightarrow \clip(\hU)\;.$$
One checks that $\ma \; 1 = 1$ and $\di |(\lab^m h)(u)| \leq (\ma^m |h|)(u)$ for all $u\in \hU$, 
$h\in \clip(\hU)$ and $m \geq 0$. It is also easy to check that $L_{f^{(0)}}^*\nu = \nu$, i.e.  $\int L_{f^{(0)}} H \, d\nu = \int H\, d\nu$ for any $H \in \clip(\hU)$.

The hyperbolicity of the flow on $\mt$ implies the existence of
 constants $c_0 \in (0,1]$ and $\gamma_1 > \gamma > 1$ such that
\begin{equation}
 c_0 \gamma^m\; d (u_1,u_2) \leq 
d (\sigma^m(u_1), \sigma^m(u_2)) \leq \frac{\gamma_1^m}{c_0} d (u_1,u_2)
\end{equation}
whenever $\sigma^j(u_1)$ and $\sigma^j(u_2)$ belong to the same  $U_{i_j}$ 
for all $j = 0,1 \ldots,m$.

Set  $\ttau = \max \{ \, \|\tau \|_0 \, , \, \Lip(\tau_{|\hU}) \, \}\;.$
Assuming that the constant $a'_0 > 0$ is sufficiently small, there exists $T = T(a_0')$ such that
\begin{equation}
T \geq \max \{ \, \|\fa \|_0 \, , \, \Lip(\fa_{|\hU}) \, , \,  \ttau\, \}\;
\end{equation}
for all $|a| \leq a'_0$. Fix $a_0' > 0$ and $T > 0$ and with these properties.
Taking the constant $C_0 > 0$ sufficiently large, we have
 $|\fa - f^{(0)}| \leq C_0|a|$ on $\hU$ for $|a| \leq a'_0$.
From now on {\bf we will assume that $a'_0$, $c_0$, $C_0$, $T$, $\gamma$ and $\gamma_1$ are fixed constants} 
with the above properties.

\section{Some properties of cylinders }
\setcounter{equation}{0}

Let again $\rr = \{ R_i\}_{i=1}^k$ be a fixed Markov family as in section 2. Define the matrix 
$A = (A_{ij})_{i,j=1}^k$  by 
$A_{ij} = 1$ if $\pp(\Int(R_i)) \cap \Int(R_j) \neq  \e$ and $A_{ij} = 0$ otherwise.  
According to \cite{kn:BR} (see section 2 there), we may assume that $\rr$ is chosen in such a way that 
$A^{M_0} > 0$ (all entries of the $M_0$-fold product of $A$  by itself are positive) for some  integer 
$M_0 > 0$. In what follows  we assume that the matrix $A$ has this property.

Given a finite string $\ii = (i_0,i_1, \ldots,i_m)$  of integers $i_j \in \{ 1, \ldots,k\}$, we will say 
that $\ii$ is {\it admissible}
if  for any $j = 0,1, \ldots,m-1$ we have $A_{i_j i_{j+1}} = 1$.  Given an admissible string $\ii$, 
denote by  $\co[\ii]$ the set
of those $x\in U$ so that $\sigma^j(x) \in \Intu(U_{i_j})$ for all $j = 0,1, \ldots,m$. The set 
$\di C[\ii] = \overline{\co[\ii]} \subset \mt$
will be called a {\it cylinder} of length $m$ in $U$, while $\co[\ii]$ will be called an {\it open cylinder} 
of length $m$. It follows from the properties of the Markov family that $\co[\ii]$ is an open dense 
subset of $C[\ii]$. Any cylinder of the form $C[i_0, i_1, \ldots, i_m, i_{m+1}, \ldots, i_{m+q}]$ 
will be called a {\it subcylinder} of $C[\ii]$ of {\it co-length} $q$.

In what follows the cylinders considered are always defined by finite admissible strings. 

The  $\sigma$-invariant probability measure $\nu$ on $U$ defined in section 2  is  a {\it Gibbs
measure} related to $g = f - P\tau$ (cf. \cite{kn:Si}, \cite{kn:R2} or \cite{kn:PP}). It follows 
from $\mbox{\rm Pr}_{\sigma}(g) = 0$,
and the properties of Gibbs measures that  there exist constants $c_2 > c_1 > 0$ such that
\begin{equation}
c_1 \leq \frac{\nu (C[\ii])}{e^{g_m(y)}} \leq c_2
\end{equation}
for any  cylinder $C[\ii]$ of length $m$ in $U$ and any $y\in C[\ii]$.
It is well-known (see \cite{kn:B}) that the core $\hU$ of  $U$ (see section 2) is a residual 
subset of $U$ and
$\nu(\hU) = 1$. Notice that for any cylinder $C[\ii]$ the set $\hC[\ii] = C[\ii] \cap \hU$ is
dense in $C[\ii]$ and $\nu(\hC[\ii]) = \nu(C[\ii])$. 

Given $x\in U_i$ for some $i$ and $r > 0$ we will denote by $B_U(x,r)$ the set of all $y\in U_i$ with 
$d (x,y) < r$.

The proof of the following proposition is straightforward.

\bs

\noindent
{\bf Proposition 3.1.} {\it If for some integer $m \geq 1$ the map 
$\sigma^m : \cc \longrightarrow \cc'$ defines a 
homeomorphism between two open cylinders $\cc$ and $\cc'$ and $w: \cc' \longrightarrow \cc$ 
is its inverse map, then 
$w(\cc'')$ is an open  subcylinder of $\cc$ of co-length $q \geq 1$ for every open subcylinder 
$\cc''$ of $\cc'$ of co-length $q$.}

\bs

Recall the constants $c_0 \in (0,1)$ and  $\gamma_1 > \gamma > 1$ from section 2,  
and fix an integer  $p_1 \geq 1$ with
\be
\rho_0 = \frac{1}{c_0\gamma^{p_1}} < \min \left\{ \frac{\diam(U_i)}{\diam (U_j)}   
: i,j = 1, \ldots, k \right\}\;.
\ee
Then clearly $\rho_0 < 1$. Set $\rho_1 = \rho_0^{1/p_1}$ and  
fix a constant $r_0 > 0$ with $2r_0 < \min \{ \diam(U_i) : i = 1, \ldots,k\}$ and for each 
$i = 1, \ldots,k$ a  point $\hz_i \in \hU_i$ such that $B_U (\hz_i, r_0) \subset \Intu(U_i)$.

The following is an immediate consequence of (2.1).

\bs

\noindent
{\bf Lemma 3.2.} {\it There exists a global constant $C_1 > 0$ such that for any cylinder 
$C[\ii]$ of length $m$ we have
$\diam (C[\ii] ) \leq C_1\, \rho_1^m$ and $\diam (C[\ii] ) \geq \frac{c_0 r_0}{\gamma_1^m}$.}
\endofproof

\bs

From now on we will assume that the local stable holonomy maps through $\mt$ are uniformly Lipschitz.
Then there exists a constant $K' > 0$ such that $d(\hh_x^y(z), \hh_x^y(z')) \leq K'\, d(z,z')$ for all 
$x,y\in \mt$ with $d(x,y) < \ep_1$ and $z,z'\in \mt\cap W^u_{\ep_1}(x)$.

Given $i = 1, \ldots, k$, according to the choice of the Markov family $\{R_i\}$,
the projection\\ $\pr_{R_i} : W_i = \phi_{[-\chi, \chi]}(R_i) \longrightarrow R_i$ along the flow $\phi_t$
is well-defined and Lipschitz.  Since the  projection $\pi_i : R_i \longrightarrow U_i$ along stable leaves 
is Lipschitz, the map $\psi_i = \pi_i\circ \pr_{R_i} : W_i  \longrightarrow W^u_{R_i}(z_i)$ is also 
Lipschitz. Thus, we may assume the constant $K' > 0$ is chosen sufficiently large so that 
$d (\psi_i(u),\psi_i (v)) \leq K'\, d (u,v)$ for all $u,v\in W_i$
and all $i = 1, \ldots,k$.

The following lemma describes the main consequences of the flow having a regular distortion along
unstable manifolds  that will be used in sections 4 and 5.

\bs

\noindent
{\bf Proposition 3.3.} {\it Assume that $\phi_t$ has a regular distortion along unstable manifolds over 
the basic set $\mt$ and that the local stable holonomy maps through $\mt$ are uniformly Lipschitz.
Then there exist global constants $0 <\rho  < 1$ and $C_1 > 0$ and a positive integer $p_0 \geq 1$ 
such that:}

(a) {\it  For  any cylinder $C[\ii] = C[i_0, \ldots,i_m]$ and any subcylinder
$C[\ii'] = C[i_0,i_{1}, \ldots, i_{m+1}]$ of $C[\ii]$ of co-length $1$  we have}
$\rho \; \diam ( C [\ii] ) \leq  \diam ( C [\ii'] ) \;.$

\ms

(b) {\it For  any cylinder $C[\ii] = C[i_0, \ldots,i_m]$ and any subcylinder
$C[\ii'] = C[i_0,i_{1}, \ldots, i_{m+1}, \ldots, i_{m+p_0}]$ of $C[\ii]$ of co-length $p_0$
we have}
$\diam(C[\ii'] ) \leq \rho\, \diam (C[\ii]) \;.$

\medskip

\noindent
{\it Proof of Proposition} 3.3. Notice that the properties of the Markov
family (and the fact that $\ii$ is admissible) imply $\sigma^m (\hC[i_0, \ldots,i_m]) = \hU_{i_{m}}$.

\medskip

(a) Set $z = \hz_{i_{m+1}}$ for brevity, and 
let $x\in \hC[\ii']$ be the point such that $\sigma^{m+1} (x) = z$.  Set $r'_0 = c_0r_0/\gamma_1 > 0$. 
Let $R = R(r'_0/K', \ep_0) > 0$ be the constant from the definition of regular distortion 
along unstable manifolds with $\delta = r'_0/K'$ and $\ep = \ep_0$ in (1.1).
Since $B_U(z,r_0) \subset \Intu (U_{i_{m+1}})$, it follows from the properties of $\sigma$ that
$B_U(\sigma^m(x), r'_0) \subset \sigma^{-1} (U_{i_{m+1}})$. Thus, for $T = \tau_m(x)$ this implies\footnote{Since
$r'_0/K' < r_0$, for any $y\in \mt \cap B^u_T(x,r'_0/K')$ one derives that 
$\sigma^j(y) \in U_{i_j}$ for all $ j = 0,1,\ldots,m+1$.} $\mt \cap B^u_T(x,r'_0/K') \subset C[\ii']$, so
$$\diam (C[\ii'] ) \geq \diam (\mt \cap B^u_T(x,r'_0/K')) \geq \frac{1}{R}\, 
\diam (\mt \cap B^u_T(x, \ep_0))\;.$$
On the other hand, 
$C[\ii] \subset \mt \cap B^u_T(x, \ep_0)$. Indeed, if $y \in C[\ii]$, then $\pp^m(y)\in R_{i_m}$ and 
$\pp^m(y) = \phi_t(\phi_T(y))$ for some $|t| \leq \chi < \ep_0/2$, where $\phi_T(y) \in W^u_{\ep_0}(\phi_T(x))$. 
Since $\diam(R_j) \leq \chi < \ep_0/2$ and $\phi_T(x) = \pp^m(x)$, we get
$d(\phi_T(x), \phi_T(y)) \leq d(\pp^m(x), \pp^m(y)) 
+ d(\pp^m(y), \phi_T(y)) \leq 2\chi < \ep_0.$
Thus,  setting  $\rho = 1/R$, we get $\rho\, \diam (C[\ii]) \leq \diam (C[\ii'])$.

\medskip

(b) Choose an arbitrary $\rho\in (0,1)$ (e.g. take the one from part (a) above), and set
$\ep = r'_0/K'$. It follows from the condition (b) in the definition of regular distortion along unstable manifolds 
that there exists $\delta \in (0,r'_0/K']$ such that
$\diam( \mt \cap B^u_T(x,\delta)) \leq \rho\; \diam( \mt \cap B^u_T(x,r'_0/K'))$
for $x\in\mt$ and $T\geq 0$.
Choose the integer $p_0 \geq 1$ so that $C_1\, K'\, \rho_1^{p_0} < \delta$.

Let $C[\ii] = C[i_0, \ldots,i_m]$ be an arbitrary cylinder and let
$C[\ii'] = C[i_0,i_{1}, \ldots, i_{m+1}, \ldots, i_{m+p_0}]$ be a subcylinder of $C[\ii]$ of co-length $p_0$.
Let $x\in \hC[\ii'] $ be the point such that $\sigma^{m+p_0} (x) = \hz_{i_{m+p_0}}$ and let $T = \tau_m(x)$.
For the cylinder $\cc' =  C[i_m, i_{m+1}, \ldots, i_{m+p_0}] \subset U_{i_m}$,
it follows from Lemma 3.2 that $\diam(\cc') \leq C_1\, \rho_1^{p_0}$, so by the choice of $p_0$,
$\diam(\cc') < \delta/K'$ and therefore $\cc' \subset B_U(\sigma^m(x), \delta/K')$. 
Next, we have $C[\ii'] \subset \mt \cap B^u_T(x, \delta)$. Indeed, if $y \in C[\ii']$, then 
$\sigma^m(y) \in \cc'$, so
$d(\sigma^m(x), \sigma^m(y)) < \delta/K'$. For  $T = \tau_m(x)$ we have 
$\phi_T(y) \in W^u_{\ep_0}(\phi_T(x))$,
so $d(\phi_T(x), \phi_T(y)) = d(\hh_{\sigma^m(x)}^{\phi_T(x)}(\sigma^m(x)),
\hh_{\sigma^m(x)}^{\phi_T(x)}(\sigma^m(y))) \leq K'\, d(\sigma^m(x), \sigma^m(y)) < \delta$.  
Thus, $C[\ii'] \subset \mt \cap B^u_T(x, \delta)$ and therefore
$\diam (C[\ii']) \leq \diam (\mt \cap B^u_T(x, \delta) ) \leq \rho \;\diam( \mt\cap B^u_T(x, r_0/K')).$
On the other hand, $B_U(\sigma^m(x), r_0) \subset \Intu (U_{i_m})$ gives
$\mt\cap B^u_T(x, r_0/K') \subset C[\ii]$, so $\diam (C[\ii']) \leq \rho\; \diam(C[\ii])$.
\endofproof

\def\Ulo{U^{(\ell_0)}}

\section{The temporal distance function}
\setcounter{equation}{0}

Throughout we assume that $\phi_t$ is a $C^2$ Axiom A flow on $M$ and $\mt$ is a basic set for 
$\phi_t$ satisfying  the condition (LNIC) stated in section 2 and such that the local holonomy maps 
along stable laminations through $\mt$ are uniformly Lipschitz.

{\bf Fix  an arbitrary point $z_0 \in \mt$ and constants $\ep_0 > 0$ and $\theta_0 \in (0,1)$ 
 with the properties described  in} (LNIC). Without loss of generality we will assume that 
$z_0 \in \Intu (U_1)$, $U_1 \subset \mt \cap W^u_{\ep_0}(z_0)$ and 
$S_1 \subset \mt \cap W^s_{\ep_0}(z_0)$, where as in section 2, $R_i = [U_i  , S_i]$ are the members
of the Markov family $\rr = \{ R_i\}_{i=1}^k$. Fix an arbitrary constant $\theta_1$ such that
$$0 <  \theta_0  < \theta_1 < 1 \;. $$

Next, fix an arbitrary orthonormal basis $e_1, \ldots, e_{n}$ in $E^u (z_0)$ and a $C^1$
parametrization $r(s) = \exp^u_{z_0}(s)$, $s\in V'_0$, of a small neighbourhood $W_0$ of $z_0$ in 
$W^u_{\ep_0} (z_0)$ such that $V'_0$ is a convex compact neighbourhood of $0$ in 
$\R^{n} \approx \mbox{\rm span}(e_1, \ldots,e_n) = E^u(z_0)$. Then $r(0) = z_0$ and
$\frac{\partial}{\partial s_i} r(s)_{|s=0} = e_i$ for all $i = 1, \ldots,n$.  Set 
$U'_0 = W_0\cap \mt $.
Shrinking $W_0$ (and therefore $V'_0$ as well)
if necessary, we may assume that $\overline{U'_0} \subset \Intu (U_1)$ and
$\left|\left\la \frac{\partial r}{\partial s_i} (s) ,  \frac{\partial r}{\partial s_j} (s) \right\ra
- \delta_{ij}\right| $ 
is uniformly small for all $i, j = 1, \ldots, n$ and $s\in V'_0$, so that
\be
\frac{1}{2} \la \xi , \eta \ra \leq \la \; d r(s)\cdot \xi \; , 
\;  d r(s)\cdot \eta\;\ra \leq 2\, \la \xi ,\eta \ra \quad, \quad
\xi , \eta \in E^u(z_0) \:,\: s\in V'_0\, ,
\ee
and
\be
\frac{1}{2}\, \|s-s'\| \leq  d ( r(s) ,  r(s')) \leq 2\, \|s-s'\|\quad, \quad s,s'\in V'_0\;.
\ee

In what follows we will construct, amongst other things, a sequence of unit vectors
$\eta_1, \eta_2, \ldots, \eta_{\ell_0}\in E^u(z_0)$. For each $\ell = 1, \ldots,\ell_0$  set 
$B_\ell = \{ \eta \in \S^{n-1} : \la \eta , \eta_\ell\ra  \geq \theta_0\}\;.$
For $t \in \R$ and $s\in E^u(z_0)$ set
$\left(I_{\eta ,t} g\right)(s) = \frac{g(s+t\, \eta) - g(s)}{t}$, $ t \neq 0\;$
({\it increment} of $g$ in the direction of $\eta$). 

\bs

\noindent
{\bf Definitions 4.1.} (a)
For a cylinder $\cc \subset U'_0$ and a unit vector $\eta \in E^u(z_0)$
we will say that a {\it separation by an $\eta$-plane occurs} in $\cc$ if there exist $u,v\in \cc$ with 
$d(u,v) \geq \frac{1}{2}\, \diam(\cc)$ such that
$ \left\la \frac{r^{-1}(v) - r^{-1}(u)}{\| r^{-1}(v) - r^{-1}(u)\|}\;,\; \eta \right\ra  \geq \theta_1\;.$

Let  $\ss_\eta$ be the {\it family of all cylinders} $\cc$ contained in $U'_0$ such that a separation by an $\eta$-plane 
occurs in $\cc$.

\ms

(b) Given an open subset $V$ of $U'_0$  which is a finite union of open cylinders and  $\delta > 0$, let
$\cc_1, \ldots, \cc_p$ ($p = p(\delta, V)\geq 1$) be the family of maximal closed cylinders in $\oV$ with
$\diam(\cc_m) \leq \delta$. For any unit vector $\eta \in E^u(z_0)$ set
$M_{\eta}^{(\delta)}(V) = \cup \{ \cc_m : \cc_m \in \ss_{\eta} \:, \: 1\leq m \leq p\}\;.$

\ms

Our aim in this section is to prove the following:

\bs

\noindent
{\bf Lemma 4.2.} ({\bf Main Lemma}) {\it There exist integers $1 \leq n_1 \leq N_0$ and $\ell_0 \geq 1$,
a sequence of unit vectors $\eta_1, \eta_2, \ldots, \eta_{\ell_0}\in E^u(z_0)$
and a non-empty open subset $U_0$ of $U'_0$ which is a finite union of open cylinders of 
length $n_1$ such that setting $\uu = \sigma^{n_1} (U_0)$ we have:}

(a) {\it For any integer $N\geq N_0$ there exist Lipschitz maps $\vl_1, \vl_2 : U \longrightarrow U$ 
($\ell = 1,\ldots, \ell_0$)  such that $\sigma^N(\vl_i(x)) = x$  for all $x\in \uu$ and $\vl_i (\uu)$ is
a finite union of open cylinders of length $N$ ($i=1,2$; $\ell = 1,2, \ldots,\ell_0$).} 

(b) {\it There exists a constant $\hd > 0$ such that for all  $\ell = 1, \ldots, \ell_0$, 
$s\in r^{-1}(U_0)$, $0 < |h| \leq \hd$ and $\eta \in B_\ell$ so that  
$s+h\, \eta \in r^{-1}(U_0\cap \mt)$ we have }
$$\left[I_{\eta,h} \left(\tau_{N}(\vl_2(\trr(\cdot ))) - \tau_{N}(\vl_1(\trr(\cdot)))\right)\right](s) 
\geq \frac{\hd}{2}\,.$$

(c) {\it We have
$\vl_i (U) \bigcap v_{i'}^{(\ell')}(U) = \e$ whenever $(i,\ell) \neq (i',\ell')$.}

(d) {\it For  any open cylinder $V$ in  $U_0$ there exists a constant 
$\delta' = \delta'(V) > 0$  such that
$$V \subset M_{\eta_1}^{(\delta)}(V) \cup M_{\eta_2}^{(\delta)}(V) \cup \ldots 
\cup M_{\eta_{\ell_0}}^{(\delta)}(V)$$ 
for all} $\delta \in (0,\delta']\;.$

\bs

Notice that if $U_0$ and $\uu$ are as in the lemma, then we must have 
$\overline{\uu} = U$.

The proof of Lemma 4.2. requires some technical preparation. We begin with a
property of the temporal distance function which follows easily from the Lipschitzness of the local stable
holonomy maps and the continuous dependence of $W^u_{\ep}(y)$ on $y \in \mt$.

\bs

\noindent
{\bf Lemma 4.3.} {\it For any $\tz\in \mt$ and any $\delta > 0$ there exists $\ep > 0$ such that 
for any $y \in W^s_{\ep}(\tz)\cap \mt$, if $y' \in W^s_{\ep}(\tz)\cap \mt$ is sufficiently close to $y$,
then for any $z,z' \in \mt \cap W^u_{\ep}(\tz)$ we have}
$$|\Delta(z',\pi_{y'}(z)) - \Delta(z', \pi_{y}(z))| < \delta\, d (z,z')\;.$$

\bs

We now proceed with the main step in the proof of the Main Lemma 4.2. This is where the
non-integrability condition (LNIC) is used.

\bs

\noindent
{\bf Lemma 4.4.} {\it There exist an integer $\ell_0 \geq 1$, 
open cylinders $\Ulo_0 \subset \ldots \subset \Uo_0$ contained in $U'_0$, and  for each 
$\ell = 1, \ldots, \ell_0$, an integer $m_\ell \geq 1$  such that the following hold:}

(i) {\it For each $\ell = 1, \ldots, \ell_0$ and each $i = 1,2$,  there exists a contracting map
$\wl_i : \Ul_0 \longrightarrow U_1$ such that $\sigma^{m_\ell}(\wl_i(x)) = x$ for all $x\in \Ul_0$,
$\wl_i : \Ul_0 \longrightarrow \wl_i(\Ul_0)$ is a homeomorphism, 
$\wl_i(\Ul_0)$ is an open cylinder in $U_1$, and 
the  sets $\overline{\wl_i(\Ul_0)}$ ($\ell = 1, \ldots, \ell_0$, $i = 1,2$)
are disjoint.}

(ii) {\it  For each $\ell = 1, \ldots, \ell_0$ there exist a number $\hd_\ell \in (0,\delta_0)$ 
and a vector $\eta_\ell\in \S^{n-1}$ such that
$| \langle  \eta_\ell , \eta_{\ell'}  \rangle |  \leq  \theta_1$ whenever $\ell' \neq \ell\;,$
and
$$\inf_{0 < |h| \leq \hd_\ell}
\left|\left[I_{\eta ,h}\left( \tau_{m_\ell}(\wl_2(r(\cdot))) 
- \tau_{m_\ell}(\wl_1(r(\cdot))) \right)\right](s)\right| \geq \hd_\ell \quad,
\quad r(s)\in \Ul_0\:, \: \eta \in B_\ell \:,$$
for all $\ell = 1, \ldots, \ell_0$, where the $\inf$ is taken over $h$ with $r(s+h \eta )\in \Ul_0$. }

(iii) {\it For  any open cylinder $V$ in  $\Ulo_0$ there exists a constant 
$\delta' = \delta'(V)\in (0,\delta_0)$  such that
$V \subset M_{\eta_1}^{(\delta)}(V) \cup M_{\eta_2}^{(\delta)}(V) \cup \ldots 
\cup M_{\eta_{\ell_0}}^{(\delta)}(V)$ for all} $\delta \in (0,\delta']\;.$

\ms

\noindent
{\bf Remark.} Notice that if $\delta' > 0$ and the open cylinder $V$ in $U'_0$ are such that
$V \subset \cup_{\ell=1}^{\ell_0} M_{\eta_\ell}^{(\delta)}(V)$
for all $\delta \in (0,\delta']$, then for any open cylinder $W \subset V$ we have 
$W \subset \cup_{\ell=1}^{\ell_0} M_{\eta_\ell}^{(\delta)}(W)$
for all $\delta \in (0,\delta'']$, where $\delta'' = \min\{ \delta', \diam(W)\}$.

\bs

\noindent
{\it Proof of Lemma} 4.4.  
Clearly if $x,y\in \S^{n-1}$ are such that $\la x,y \ra  \leq \theta_1$, 
then $\|x - y \| \geq \sqrt{2(1-\theta_1)}$. Thus there exists a positive integer $\ell'_0$, 
depending on $n$ and $\theta_1$ only, such that for any finite set 
$\{x_1, \ldots, x_k\} \subset \S^{n-1}$ with $\la x_i, x_j\ra \leq \theta_1$ for all
$i \neq j$ we have $k \leq \ell'_0$. {\bf Fix $\ell'_0$ with this property.} 

As another preparatory step, fix $2\ell'_0$ distinct points
$x_i^{(\ell)} \in \Intu( U_1) \setminus \overline{U'_0}$ 
($\ell = 1, \ldots, \ell'_0 \; ; \:  i = 1,2$) and for each $x_i^{(\ell)}$ 
{\bf fix a small open neighbourhood} ${\tVl}_i$ of $x_i^{(\ell)}$ in 
$\Intu (U_1)$ such that the sets 
$\overline{{\tVl}_i}$ ($\ell = 1, \ldots, \ell'_0\;; \: i = 1,2$) are disjoint and contained
in $\Intu (U_1)$.

We will construct  the required  objects by induction.

\ms

\noindent
{\bf Step 1.} Clearly there exists a unit vector $\eta_1 \in E^u(z_0)$ tangent to $\mt$ at $z_0$.
It then follows from the condition (LNIC) and the choice of $z_0$ that there exist $\tz = r(\ts) \in U'_0$, 
$\ty_1, \ty_2 \in W^s_{R_1}(\tz)$ (so $\ty_1, \ty_2 \in \mt$)
with $\ty_1 \neq \ty_2$, $\delta'_1 > 0$ and $\ep'_1 > 0$ such that 
\be
|\Delta( r(s + h\, \eta), \pi_{\ty_1}(r(s)))
- \Delta( r(s + h\, \eta), \pi_{\ty_2}(r(s)))| \geq \delta'_1\, | h|
\ee
for all $r(s)\in U'_0$ with $\dist(\tz,r(s)) < \ep'_1$, $\eta \in B_1$ and $h\in \R$ with $|h|< \ep'_1$ and
$r(s+ h\, \eta) \in U'_0$. 
We will assume that $\ep'_1 > 0$ is so small that $B_U(\tz,\ep'_1) \subset U'_0$.

Since $\tVo_1$ and $\tVo_2$ are open subsets of $U$ having common points with 
$\Lambda$, it follows that $\pp^{m}(\tVo_1)$ and $\pp^{m}(\tVo_2)$  fill in $R_1$ densily 
as $m \to \infty$.  Using this,  it follows that  taking $m_1 \geq 1$ large enough we can find 
$y'_1\in W^s_{R_1}(\tz) \cap \pp^{m_1}(\tVo_1)$  arbitrarily close to $\ty_1$ and 
$y_2' \in  W^s_{R_1}(\tz) \cap \pp^{m_1}(\tVo_2)$ arbitrarily close to
$\ty_2$.  By Lemma 4.3
we can make this choice so  that for any $i =1,2$ and any $z = r(s)$, $h$ and $\eta$ as above, we have
\begin{eqnarray*}
| \Delta (r(s + h\, \eta), \pi_{y'_i}(r(s))) - \Delta (r(s + h\, \eta), \pi_{\ty_i}(r(s)))| \leq   \frac{\delta'_1\, |h|}{4}.
\end{eqnarray*}
Combining this with (4.3) one gets
\begin{eqnarray*}
&         & |\Delta (r(s + h\, \eta), \pi_{y'_1}(r(s))) - \Delta (r(s + h\, \eta), \pi_{y'_2}(r(s)))|\\
& \geq & |\Delta (r(s + h\, \eta), \pi_{\ty_1}(r(s)))|  - |\Delta (r(s + h\, \eta), \pi_{\ty_2}(r(s)))|\\
&      & - |\Delta (r(s + h\, \eta), \pi_{\ty_1}(r(s)))|  - |\Delta (r(s + h\, \eta), \pi_{y'_1}(r(s)))|\\
&      & - | \Delta (r(s + h\, \eta), \pi_{y'_2}(r(s))) - \Delta (r(s + h\, \eta), \pi_{\ty_2}(r(s)))| \geq \frac{\delta'_1\, |h|}{2}\;.
\end{eqnarray*}

Thus, there exists an open cylinder $\Uo_0$ contained in $B_U(\tz,\ep'_1) \subset U'_0$ 
with $\tz \in \Uo_0$ such that
\begin{equation} 
|\Delta (r(s + h\, \eta), \pi_{y'_1}(r(s))) - \Delta (r(s + h\, \eta), \pi_{y''_1}(r(s)))| 
\geq \hd_1\, |h| \;
\end{equation}
whenever $r(s) \in \Uo_0$, $ \eta \in B_1$, $| h| \leq \hd_1$ and $r(s+h \eta) \in \Uo_0$,
where $\hd_1 = \min\{ \delta'_1/2, \ep'_1\}$. {\bf Fix $m_1 \geq 1$,  $y'_1$ and $y'_2$ 
with these properties}.

Let $\ooo_1$ be a small open neighbourhood of $y'_1$ in  $W^u_{R_1}(y'_1) \cap \pp^{m_1}(\tVo_1)$
and let\\ $\fo_1: \ooo_1 \longrightarrow \fo_1(\ooo_1) \subset \tVo_1$ be a 
Lipschitz homeomorphism (local inverse of $\pp^{m_1}$) such that  $\pp^{m_1}(\fo_1(z)) = z$  
for all $z\in \ooo_1$. 
Shrinking $\Uo_0$ if necessary, we may assume that $\pi_{y'_1} (\Uo_0) \subset \ooo_1$. 
Now define a Lipschitz homeomorphism
$\wo_1 : \Uo_0 \longrightarrow \wo_1 (\Uo_0)\subset  \tVo_1$ by  $\wo_1(x) = \fo_1(\pi_{y'_1}(x))$.
We then have $\pp^{m_1}(\wo_1(x)) = \pi_{y'_1} (x)$ and therefore
$\sigma^{m_1}(\wo_1(x)) = x$ for all $x\in \Uo_0$.  Moreover, 
$\Lip(\wo_1) \leq \frac{1}{c_0\, \gamma^{m_1}}$, so assuming $m_1$ is sufficiently large,
$\wo_1$ is contracting and $\wo_1(\Uo_0)$ is a cylinder.

In the same way  one  constructs a Lipschitz homeomorphism 
$\wo_2 : \Uo_0 \longrightarrow \wo_2 (\Uo_0) \subset\tVo_2$ 
(replacing $\Uo_0$ by a smaller open cylinder if necessary; by Proposition 3.1,  
$\wo_1(\Uo_0)$ will continue 
to be a cylinder)  such that $\wo_2(\Uo_0)$ is a cylinder and 
$\pp^{m_1}(\wo_2(x)) = \pi_{y'_2}(x)$ for all $x\in \Uo_0$.  
Then $\sigma^{m_1}(\wo_2(x)) = x$ for all $x\in \Uo_0$.  

Now  for $z = r(s) \in \Uo_0$, $\eta\in B_1$ and $h\in \R$ with $r(s + h\, \eta) \in \Uo_0$ we get
\begin{eqnarray*}
&    & [\tau_{m_1}(\wo_2(r(s+ h\, \eta))) - \tau_{m_1}(\wo_1(r(s+ h\, \eta)))]
  - [\tau_{m_1}(\wo_2(r(s))) - \tau_{m_1}(\wo_1(r(s)))]\\
& = & [\tau_{m_1}(\wo_1(r(s))) - \tau_{m_1}(\wo_1(r(s+ h\, \eta)))]
- [\tau_{m_1}(\wo_2(r(s)))  - \tau_{m_1}(\wo_2(r(s+ h\, \eta)))]\\
& = & \Delta ( \pp^{m_1}(\wo_1(r(s+ h\, a))), \pp^{m_1}(\wo_1(r(s)) )
  - \Delta( \pp^{m_1}(\wo_2(r(s + h\, \eta))) , \pp^{m_1}(\wo_2(r(s)) )\\
& = & \Delta ( \pi_{y'_1}(r(s+ h\, \eta)), \pi_{y'_1}(r(s)) ) 
- \Delta (\pi_{y'_2}(r(s + h\, \eta) ),  \pi_{y'_2}(r(s) )) \\
& = &  \Delta ( r(s+ h\, \eta), \pi_{y'_1}(r(s)) ) - \Delta ( r(s+ h\, \eta), \pi_{y'_2}(r(s)) )\;.  
\end{eqnarray*}
This and  (4.4) give 
$| [I_{\eta,h} (\tau_{m_1}(\wo_2(r(\cdot))) - \tau_{m_1}(\wo_1(r(\cdot))))](s) | \geq \hd_1$
whenever $r(s)\in \Uo_0$, $\eta \in B_1$, $0 < |h| \leq \ep'_1$ and $r(s+h\, \eta) \in \Uo_0$.

In this way we have completed the first step in our recursive construction. Whether we 
need to make  more steps or not depends on which of the following two alternatives takes place.

\ms

\noindent
{\bf Alternative 1.A.} There exist an  open cylinder $V$ contained in $\Uo_0$ and  a constant
$\delta_1 \in (0,\delta_0)$  such that $M_{\eta_1}^{(\delta)}(V) \supset V$  for all $ \delta\in (0, \delta_1]\;.$

\ms

\noindent
{\bf Alternative 1.B.} Alternative 1.A does not hold.

\ms

In the case of Alternative 1.A we simply  terminate the recursive construction at this stage replacing
 $\Uo_0$ by $V$.

If Alternative 1.B takes place, we need to make at least one more step.

\ms

\noindent
{\bf Inductive Step.} Suppose that for some $j \geq 1$ we have constructed 
open cylinders $U^{(j)}_0 \subset \ldots \subset \Uo_0$ contained in $U'_0$,
and for each $\ell = 1, \ldots, j$, an integer $m_\ell \geq 1$  and a vector
$\eta_\ell \in \sn$ such that the conditions (i) and (ii) in the lemma are fulfilled with 
$\ell_0$ replaced by $j$.

There are two alternatives again.

\ms

\noindent
{\bf Alternative j.A.} There exist an open cylinder $V$ contained in $\Uj_0$ and
$\delta_j \in (0,\delta_{j-1}]$  such that
$$M_{\eta_1}^{(\delta)}(V) \cup M_{\eta_2}^{(\delta)}(V)  \cup \ldots \cup M_{\eta_j}^{(\delta)}(V) 
\supset V  \quad, \quad \delta\in (0, \delta_j]\;.$$

\ms

\noindent
{\bf Alternative j.B.} Alternative j.A  does not hold.

\ms

In the case of Alternative j.A we  terminate the recursive construction at this stage replacing
 $\Uj_0$ by $V$.

Next, assume that Alternative j.B takes place. One then needs to complete

\ms

\noindent
{\bf Step j+1.} Construct an  open cylinder $U^{(j+1)}_0$ contained in  $\Uj_0$,
an integer $m_{j+1} \geq 1$ and a unit vector $\eta_{j+1}\in E^u(z_0)$  such that
the conditions (i) and (ii) in the lemma are fulfilled with $\ell_0$ replaced by $j+1$.

Given an open cylinder $V$ in $\Uj_0$ and $\delta > 0$, set 
$A_\delta (V) = M_{\eta_1}^{(\delta)}(V) \cup M_{\eta_2}^{(\delta)}(V) \cup  \ldots 
\cup M_{\eta_j}^{(\delta)}(V)$.
It follows from Alternative j.B that for any open cylinder $V$ contained in $\Uj_0$
we have $V \setminus A_\delta(V) \neq \e$ for arbitrarily small $\delta$.

Notice that, since there are only countably many cylinders in $\Uj_0$, there exists a decreasing sequence
$\delta_j = \mu_0 > \mu_1 > \mu_2 > \ldots > \mu_k > \ldots$ converging to zero such that for any cylinder
$\cc$ in $\Uj_0$ we have $\diam(\cc) = \mu_k$ for some $k$. Then for any open cylinder $V\subset \Uj_0$ and any
$m \geq 1$ we have $A_\delta(V) = A_{\mu_m}(V)$ whenever $\mu_m \leq \delta < \mu_{m-1}$.

Let $V = \Uj_0$ and for any $m \geq 1$ consider the compact subset
$F_m = \bigcap_{k \geq m} A_{\mu_k}(V)\;$
of $\overline{V}$. Clearly $F_{m'} \subset F_m$ whenever $m' < m$.
We claim that $\Intu (F_m) = \e$ for all $m\geq 1$. Indeed, assume that $\Intu (F_m) \neq \e$
for some $m' \geq 1$; then there exists a non-empty open cylinder $W$ contained in $F_{m'}$. There exists
$k' \geq 1$ with $\mu_{k'} = \diam(W)$. Setting $m = \max\{ m', k' \} + 1$, we have $W \subset F_m$.
Moreover, for $0 < \delta \leq \mu_{m}$ we have 
\be
M_{\eta_i}^{(\delta)}(V) \cap W \subset M_{\eta_i}^{(\delta)}(W)\quad, \quad i = 1,\ldots, j\;.
\ee
Consequently, for all $0 < \delta \leq \mu_{m}$ we have $A_\delta(V) \cap W \subset A_\delta(W)$. 
Now $W \subset F_{m}$  implies
$W = F_m \cap W = \bigcap_{k \geq m} A_{\mu_k}(V) \cap W
\subset \bigcap_{k \geq m} A_{\mu_k}(W)\;.$
Thus, $W \subset A_{\mu_k}(W)$ for all $k \geq m$, which means that $W \subset A_{\delta}(W)$
for all $0 < \delta \leq \mu_{m}$. This is a contradiction with Alternative j.B.
 
Hence $\Intu (F_m) = \e$ for all $m\geq 1$. Thus, $\Uj_0 \setminus F_m$ are open and dense subsets 
of $\Uj_0$, so $G = \cap_{m=1}^\infty (\Uj_0 \setminus F_m)$ is a residual (even a $G_\delta$) 
subset of $\Uj_0$. The properties of $\hU$ now imply $G \cap \hU \neq \e$.  

Choose an arbitrary $\hz \in G\cap \hU$ and fix it. Given any $m\geq 1$, we have $\hz \notin F_m$, so
there exists $\mu'_m \in (0,\mu_m]$ with $\hz \notin A_{\mu'_m}(V)$, i.e.
 $\hz \notin \cup_{\ell=1}^j  M_{\eta_\ell}^{(\mu'_m)}(V)$. So, if $\ccm$ is the maximal cylinder in 
$ \overline{V}$  with $\diam(\ccm ) \leq \mu'_m$ such that $\hz \in \ccm$, then
$\ccm\notin \ss_{\eta_\ell}$ for any $\ell = 1, \ldots, j$.

Fix $m \geq 1$  for a moment, and let $u_m, v_m \in\ccm$ be  such that $d(u_m, v_m) = \diam(\ccm)$.
Since $\hz\in \ccm$, we may assume $d(u_m,\hz) \geq \frac{1}{2}\, \diam(\ccm)$.
Then $\left\la \frac{r^{-1}(u_m) - r^{-1}(\hz)}{\| r^{-1}(u_m) - r^{-1}(\hz)\|}\, , \, 
\eta_\ell  \right\ra  < \theta_1\;$
for all $\ell = 1, \ldots,j$, since $\ccm\notin \ss_{\eta_\ell}$.
Clearly, $u_m \to \hz$ as $m\to \infty$. Choose a subsequence $\{u_{m_p}\}$ so that
$\eta_{j+1} = \lim_{p\to\infty} \frac{r^{-1}(u_{m_p}) - r^{-1}(\hz)}{\|r^{-1}(u_{m_p}) 
- r^{-1}(\hz)\|}  \in \S^{n-1}$
exists.  Let $\hz = r(\hs)$; then $\eta'_{j+1} = dr(\hs)\cdot \eta_{j+1}$ is tangent to
$\mt$ at $\hz$, and according to the above, $\left\la \eta_{j+1}, \eta_\ell  \right\ra  \leq \theta_1\;$
for all $\ell = 1, \ldots, j$.

Repeating the argument from the proof of (4.4) in Step 1, one derives that there exist 
$\tz = r(\ts) \in V  = \Uj_0$, $\ty_1, \ty_2 \in W^s_{R_1}(\tz) \setminus U_1$ with $\ty_1 \neq \ty_2$,
$\delta'_{j+1} \in (0,\hd_j)$  and $\ep'_{j+1} > 0$ such that 
$$|\Delta( r(s + h\, \eta), \pi_{\ty_1}(r(s)))
- \Delta( r(s + h\, \eta), \pi_{\ty_2}(r(s)))| \geq \delta'_{j+1}\, | h|$$
for all $r(s)\in V $ with $d (\tz,r(s)) < \ep'_{j+1}$, $\eta \in B_{j+1}$ and $h\in \R$ with 
$|h|< \ep'_{j+1}$ and $r(s+h \eta) \in \Uj_0$. 
We will assume that $\ep'_{j+1} > 0$ is so small that $B_U(\tz,\ep'_{j+1})  \subset V = \Uj_0$.

Then, again as in Step 1, one constructs an open cylinder $\Ujj_0 \subset \Uj_0$ and
contracting homeomorphisms $\wjj_i : \Ujj_0 \longrightarrow \wjj_i (\Ujj_0) \subset \tVjj_i$ ($i = 1,2$)
with $\sigma^{m_{j+1}}(\wjj_i(x)) = x$ for all  $x\in \Ujj_0$ and such that
$| [I_{\eta ,h}(\tau_{m_{j+1}}(\wjj_2(r(\cdot))) - \tau_{m_{j+1}}(\wjj_1(r(\cdot))))](s) |\geq \hd_{j+1}$
for all $r(s)\in \Ujj_0$, $ \eta \in B_{j+1}$, $0 < | h| \leq \hd_{j+1}$ with $r(s+h \eta) \in \Ujj_0$.
This completes Step $j+1$.

By the definition of $\ell'_0$, it is clear that this inductive procedure terminates after not 
more than $\ell'_0$ steps. That is, for some $\ell_0 \leq j'_0$ the Alternative $\ell_0.A$
holds, and then we  terminate the construction at that step. \endofproof

\bs

In what follows we use the objects constructed in Lemma 4.4.
Set 
$\di \hd = \min_{1\leq \ell\leq \ell_0} \hd_j$, 
$\di n_0 = \max_{1\leq \ell\leq \ell_0} m_\ell$, 
and fix an arbitrary point $\hz_0 \in U_0^{(\ell_0)}\cap \hU$.

\medskip

\noindent
{\bf Lemma 4.5.} {\it There exist an integer $n_1 \geq 1$ and 
an open neighbourhood $U_0$ of $\hz_0$ in $\Ulo_0$ 
such that $\Intu (U) = \sigma^{n_1}(U_0)$,
$\sigma^{n_1} : U_0 \longrightarrow  \sigma^{n_1}(U_0)$ is a homeomorphism
and $U_0$ is a finite union of open cylinders of length $n_1$.}

\medskip

\noindent
{\it Proof of Lemma 4.5.}  Let $U_0^{(\ell_0)} = \co[\ii] = \co[i_0, \ldots,i_m]$. By construction,
$\hz_0 \in \hC [\ii] \subset \Ulo_0$. For the matrix $A$ we have $A^{M_0} > 0$
for some integer $M_0 \geq 1$ (see the beginning of section 3), so for each
$j = 1, \ldots,k$ there exists an admissible  string
$s^{(j)} = ( i_{m+1}^{(j)}, \ldots, i_{m+M_0}^{(j)} )$
such that $ i_{m+M_0}^{(j)} = j$ and $A_{i_m  i_{m+1}^{(j)}} = 1$. Fix an arbitrary string 
$s^{(j)}$ with this property. 
For the particular  $j'$ such that $\sigma^{m+M_0}(\hz_0) \in U_{j'}$, choose
$s^{(j')}$ in such a way that $\hz_0 \in C[\ii; s^{(j')}]$.

Now set $n_1 = m + M_0$ and $U_0 = \cup_{j=1}^k \co [\ii; s^{(j)}]\;.$
Clearly, $U_0$ is an open subset of $\Ulo_0$ containing $\hz_0$
and one shows easily that $\sigma^{n_1} : U_0 \longrightarrow \sigma^{n_1} (U_0)$ is a homeomorphism
and $\Intu (U) =  \bigcup_{j=1}^k \sigma^{n_1}(\co [\ii; s^{(j)}]) = \sigma^{n_1}(U_0)$.
 \endofproof

\bs

Using the above lemma, {\bf fix  $n_1 > 0$ and an open neighbourhood $U_0$ of $\hz_0$ in 
$U_0^{(\ell_0)}$} such that
$U_0$ is a finite union of open cylinders of length $n_1$, $\uu = \sigma^{n_1}(U_0) = \Intu (U )$
and $\sigma^{n_1} : U_0 \longrightarrow  \uu$ is a homeomorphism.
The inverse homeomorphism $\psi : \uu \longrightarrow U_0$ is Lipschitz, so it has a Lipschitz extension
\begin{equation}
\psi : U \longrightarrow \overline{U_0}\:\: \mbox{\rm  such that }\:\: 
\sigma^{n_1}(\psi(x)) = x \:\:, \:\:x\in \uu\;.
\end{equation}
Then $\trr(s) = \sigma^{n_1}(r(s))$,  $s\in V_0\;,$ where $V_0 = r^{-1}(U_0) \subset V'_0$,
gives a Lipschitz parametrization of $\uu$ with $\psi(\trr(s)) = r(s)$ for all $s \in V_0$.
Finally, set
\be
\Vl_i = \wl_i(U_0)  \subset \tVl_i\quad, \quad i = 1,2\; ; \: \ell = 1, \ldots,\ell_0\;.
\ee
It follows from the choice of $U_0$, the properties of $\wl_i$ (see (i) in Lemma 4.4) and Proposition 3.1
that $\Vl_i$ is a finite union of open cylinders of lengths $n_1 + m_\ell$.

The following  two lemmas are proved essentially by using arguments from \cite{kn:D2} and Lemma 4.4 above.
We omit most of the details.

\bigskip

\noindent
{\bf Lemma 4.6.}   
{\it For every $\delta '' > 0$ there exists an integer $n_2 > 0$ such that for any $m  \geq n_0+ n_2$, 
any $\ell = 1, \ldots, \ell_0$ and $i = 1,2$ there exist contracting maps  
$\tvl_i : \Vl_i  \longrightarrow U $  with $\sigma^{m-m_\ell}(\tvl_i(w)) = w$ for all 
$w\in \Vl_i $ such that
 \be
 \Lip (\tau_{m-m_\ell} \circ \tvl_i )\leq \delta'' \:\:\: \mbox{\rm   on }\:\:\:  \Vl_i \;,
 \ee
 $\tvl_i (\Vl_i )$ is a finite union of open cylinders of length $n_1 + m$ and
 $\overline{\tvl_i (\Vl_i )} \bigcap \overline{\tilde{v}_{i'}^{(\ell')}(V_{i'}^{(\ell')})} = \e$ whenever
 $(i,\ell) \neq (i',\ell')$.} $\endofproof$

\bigskip

Set $\delta'' = \frac{c_0\, \hd}{8} \; ,$
{\bf fix $n_2 = n_2(\delta'') > 0$ with the properties listed  in Lemma 4.6}, and  denote
$N_0 = n_0 + n_1 + n_2\;.$

\bs


\noindent
{\it Proof of Lemma} 4.2. 
Let $N \geq N_0$.  Then $m = N - n_1 \geq n_0 + n_2$, so by Lemma 4.6 for any $\ell = 1, \ldots,\ell_0$ 
and any  $i = 1,2$ there exists a  contracting homeomorphism
\begin{equation}
\tvl_i : \Vl_i \longrightarrow \tvl_i (\Vl_i ) \subset U\:\: \mbox{\rm  with } \:\: 
\sigma^{N- m_\ell - n_1}(\tvl_i(w)) = w \:\: , \:\: w\in \Vl_i \; ,
\end{equation}
such that (4.8) holds with $m = N - n_1$ and $\delta''$ as above. Moreover, we can choose
the maps $\tvl_i$ so that 
$\overline{\tvl_i (\Vl_i )} \bigcap \overline{\tilde{v}_{i'}^{(\ell')}(V_{i'}^{(\ell')})} = \e$ 
whenever $(i,\ell) \neq (i',\ell')$. Now define Lipschitz maps
\begin{equation}
\vl_i : U  \longrightarrow U \quad \mbox{\rm such that} \quad \vl_i(x) 
= \tvl_i(\wl_i(\psi(x)))\:, \: x\in \uu\;.
\end{equation}
It follows immediately from the above that $\overline{\vl_i (U)} \cap \overline{v_{i'}^{(\ell')}(U)} = \e$ 
whenever $(i,\ell) \neq (i',\ell')$, while Proposition 3.1 shows that each $\vl_i (\uu)$ is
a finite union of  open cylinders of length $N$.

Moreover, for any $x\in \uu$, according to  (i) in Lemma 4.4, (4.10)  and (4.9), we have
$\sigma^{N-n_1}(\vl_i(x)) = \sigma^{m_\ell}( \sigma^{N-m_\ell-n_1}(\vl_i(x))) =
\sigma^{m_\ell}  (\wl_i(\psi(x))) = \psi(x)\;,$
which is the same for all $\ell$ and  $i$. Consequently,
$\sigma^p(\vl_1(x)) = \sigma^p(\vl_2(x))$ for all $p \geq N-n_1$ and $x\in \uu $. Thus,
$\tau_N(\vl_2(x)) - \tau_N(\vl_1(x)) = \tau_{N-n_1}(\vl_2(x))- \tau_{N-n_1}(\vl_1(x))$ for $x \in \uu$,
and given $\eta \in B_\ell$ and $h > 0$, we have
\begin{eqnarray*}
\left[ I_{\eta ,h}\left(\tau_N(\vl_2(\trr(\cdot))) - \tau_N(\vl_1(\trr(\cdot)))\right)\right](s)
 = \left[ I_{\eta,h}\left(\tau_{N-n_1}(\tvl_2(\wl_2(r(\cdot)))) 
 -  \tau_{N-n_1}(\tvl_1(\wl_1(r(\cdot))))\right) \right](s)\;. 
\end{eqnarray*}
Now Lemma 4.4(ii), (4.2) and (4.8) imply
$ \left| \left[I_{\eta ,h} (\tau_N(\vj_2(\tr(\cdot))) 
- \tau_N(\vj_1(\tr(\cdot))))\right](s) \, \right|\geq \frac{\hd}{2}\;.$
$\endofproof$

\def\hcc{\widehat{\cc}}

\section{Dolgopyat   operators}
\setcounter{equation}{0}

In this section we prove Theorem 1.1.

Throughout we assume that $\mt$ is a basic set for a $C^2$ Axiom A flow 
$\phi_t : M \longrightarrow M$ satisfying the condition (LNIC) which has regular distortion along 
unstable manifolds over $\mt$ and uniformly Lipschitz local stable holonomy maps.
We use the notation from section 2, in particular the fixed real-valued function 
$f \in \clip(\hU)$, the function 
$g = f - P\, \tau$, where $P\in \R$ is such that $\Pr_{\sigma} (g) = 0$, and the $\sigma$-invariant
probability measure $\nu$ on $\mt$ such that $L_{f^{(0)}}^* \nu = \nu$.

The central point here is to prove the $L^1$-contraction property of 
the normalized operator $\lab$ with respect to the Gibbs measure $\nu$ and the norm $\| h\|_{\lip,b}$.  

\bigskip

\noindent
{\bf Theorem 5.1.} {\it There exist a positive integer $N$ and
constants $\hrho \in (0,1)$ and $a_0 > 0$ such that for any $a, b\in \R$ with $|a|\leq a_0$ 
and $|b|\geq 1/a_0$ and any $h\in \clip(\hU)$ with $\|h \|_{\lip,b} \leq 1$ we have
$\di\int_{U} |\lab^{N m}h|^2\; d\nu \leq \hrho^m \; $ for every positive integer $m$.}

\bs

Theorem 1.1 is derived from the above in the same way as in \cite{kn:D2} (see also the proof of
Corollary 3.3(a) in \cite{kn:St1}). Indeed, the assumptions of Theorem 1.1 are exactly the ones we have
in this section.

\ms

Define a {\it new metric} $D$ on $\hU$ by 
$$D(x,y) = \min \{ \diam(\cc) : x,y\in \cc\:, \: \cc \: \mbox{\rm a cylinder contained in }\, U_i \}$$
if $x,y \in U_i$ for some $i = 1, \ldots,k$, and $D(x,y) = 1$ otherwise.
Recall that $\diam(U_i) < 1$ for all $i$ by the choice of the Markov family.

The proof of the following lemma is straightforward.

\bs

\noindent
{\bf Lemma 5.2.} {\it (a) $D$ is a metric on $\hU$, and 
if $x,y\in \hU_i$ for some $i$,  then $d(x,y) \leq  D(x,y)$.}

\ms

(b) {\it For any cylinder $\cc$ in $U$ the characteristic function $\chi_{\hcc}$ of
$\hcc$ on $\hU$ is Lipschitz with respect to $D$ and $ \Lip_D(\chi_\cc) \leq 1/\diam(\cc)$.}
\endofproof

\bs

We will denote by $\clip_D(\hU)$ the {\it space of all Lipschitz functions $h : \hU \longrightarrow \C$
with respect to the metric} $D$ on $\hU$ and by $\Lip_D(h)$ the {\it Lipschitz constant} of $h$ with respect to $D$.

Given $A > 0$, denote by $K_A(\hU)$  {\it the set of all functions } $h\in \clip_D(\hU)$ 
such that  $h > 0$ and
$\frac{|h(u) - h(u')|}{h(u')} \leq A\, D (u,u')$ for all $u,u' \in \hU$ that belong to the 
same $\hU_i$  for some $i = 1, \ldots,k$.
Notice that $h\in K_A(\hU)$ implies $|\ln h(u) - \ln h(v)| \leq A\; D (u,v)$ and therefore
$e^{-A\; D (u,v)} \leq \frac{h(u)}{h(v)} \leq e^{A \; D (u,v)}$ for all $u, v\in \hU_i $, $i = 1, \ldots,k\;.$

\bs

Theorem 5.1 is derived from  the following lemma which  is the analogue of Lemma $10''$ in \cite{kn:D2}.
It should be stressed that replacing the standard metric\footnote{In fact, it is not clear at all whether 
a similar lemma will be true for general basic sets on manifolds of arbitrary dimension with the metric 
$d$ in the place of $D$. }  $d$ by the metric $D$ is significant here. 

\bs

\noindent
{\bf Lemma 5.3.} {\it There exist a positive integer $N$ and constants
$\hrho = \hrho(N) \in (0,1)$, $a_0 = a_0(N) > 0$, $b_0 = b_0(N) > 0$ and $E \geq 1$
such that for every $a, b\in \R$ with $|a|\leq a_0$, $|b| \geq b_0$, there exists a
finite family $\{ \nn_J\}_{J\in \J}$ of  operators
$\nn_J  = \nn_J(a,b) : \clip_D (\hU) \longrightarrow \clip_D (\hU)\;,$ 
where $\J = \J(a,b)$ is a finite set depending on $a$ and $b$, with the following properties:}

(a) {\it The operators $\nn_J$ preserve the cone} $K_{E|b|} (\hU)$ ;

(b) {\it For all $H\in K_{E|b|}(\hU)$ and $J \in \J$ we have}
$\di \int_{ \hU} (\nn_J H )^2 \; d\nu \leq \hrho \; \int_{\hU} H^2 \; d\nu\; .$

(c) {\it If $h, H\in \clip_D (\hU)$ are such that $H\in K_{E|b|}(\hU)$, $|h(u)|\leq H(u)$ 
for all   $u \in \hU$ and \\
$| h(u) - h(u') | \leq E|b| H(u')\, D (u,u')$ whenever $u,u'\in \hU_i $ for some $i = 1, \ldots,k$,  
then there exists 
$J\in \J$  such that $|\lab^N h(u)|\leq (\nn_J H)(u)$ for all $u \in \hU$   and
$$|(\lab^N h)(u) - (\lab^N h)(u')|\leq E|b| (\nn_J H)(u')\, D (u,u')$$
whenever $u,u'\in \hU_i $ for some $i = 1, \ldots,k$.}

\bs

The remainder of this section if devoted to the proof of Lemma 5.3.
We begin with a technical lemma containing two specific versions of Lasota-Yorke type of inequalities.
Its proof is given in the Appendix.
Here we use the constants $a'_0 > 0$, $T > 0$ and $\gamma_1 > \gamma > 1$ from section 2.

\bs

\noindent
{\bf Lemma 5.4.}  {\it There exists a constant $A_0 > 0$ 
such that for all $a\in \R$ with $|a|\leq a'_0$ the following hold:}

\ms

(a)  {\it If $H \in K_B(\hU)$ for some $B > 0$, then 
$\frac{|(\ma^m H)(u) - (\ma^m H)(u')|}{(\ma^m H)(u')} \leq
A_0 \, \left[ \frac{B}{\gamma^m} + \frac{T}{\gamma-1}\right]\, D (u,u')$
for all $m \geq 1$ and all $u,u'\in U_i$, $i = 1, \ldots, k$.}

\ms

(b) {\it If the functions $h$ and  $H$ on $\hU$  and $B > 0$  are such that $H > 0$ on $\hU$ and 
$|h(v) - h(v')| \leq B\, H(v')\, D (v,v')$ for any $v,v'\in \hU_i$, $i = 1, \ldots,k$, 
then for any integer $m \geq 1$ and any $b\in \R$ with  $|b|\geq 1$ we have 
$ |\lab^m h(u) - \lab^m h(u')| \leq  A_0\,\left[ \frac{B}{\gamma^m} \, (\ma^m H)(u') + |b|\, (\ma^m |h| )(u')\right]\,
D (u,u')$
whenever $u,u'\in \hU_i$ for some $i = 1, \ldots,k$.}

\bs

As in the beginning of section 4, fix  an arbitrary point $z_0 \in \Intu(U_1)$ and constants 
$\ep_0 > 0$ and $0 < \theta_0 < \theta_1 < 1$, a $C^1$
parametrization $r(s)$ of a small neighbourhood $W_0$ of $z_0$ in 
$W^u_{\ep_0} (z_0)$ with (4.1) and (4.2). Next, fix
integers $1 \leq n_1 \leq N_0$ and $\ell_0 \geq 1$, 
unit vectors $\eta_1, \eta_2, \ldots, \eta_{\ell_0}\in E^u(z_0)$
and a non-empty open subset $U_0$ of $W_0$ which is a finite union of open cylinders of 
length $n_1$ with the properties described in Lemma 4.2.
We will also use the set $\uu = \sigma^{n_1}(U_0) = \Intu (U)$ and
the constants $\rho \in (0,1)$ and $p_0 \geq 1$ from Lemma 3.3. 
Since $\sigma^{n_1} : U_0 \longrightarrow \uu$ is one-to-one, it has an inverse
map $\psi : \uu \longrightarrow U_0$, which is Lipschitz.


We will now impose certain condition on the numbers $N$, $\ep_1$, $b$
and $\mu$ that will be used throughout. Where these conditions come from
will become clear later on.

Set $E = \max\left\{  4A_0\;, \; \frac{2 A_0\, T}{\gamma-1}\; \right\},$
where $A_0 \geq 1$ is the constant from Lemma 5.4, and {\bf fix an integer} $N \geq N_0$ such that
\be
\gamma^N \geq \max \left\{ \: 6A_0 \; , \; \frac{200\, \gamma_1^{n_1}\,A_0}{c_0^2} \; , \; 
\frac{512\, \gamma^{n_1}\, E}{c_0\, \hd\, \rho}  \; \right\} \;.
\ee
Now fix maps $\vl_i : U \longrightarrow U$ ($\ell = 1, \ldots, \ell_0$, $i = 1,2$) with
the properties (a), (b), (c) and (d) in Lemma 4.2. 
In particular, (c) gives
\be
\overline{v_{i}^{(\ell)}(U)} \cap \overline{v_{i'}^{(\ell')}(U)} = \e \quad ,\quad (i,\ell) \neq (i',\ell')\;.
\ee

Since $U_0$ is a finite union of open cylinders contained in $U^{(\ell_0)}_0$ 
it follows from Lemma 4.2(d)  that there exist a  constant
$\delta' = \delta'(U_0) \in (0,\delta_0)$  such that
\be
 M_{\eta_1}^{(\delta)}(U_0) \cup \ldots\cup M_{\eta_{\ell_0}}^{(\delta)}(U_0)
\supset U_0  \quad , \quad \delta\in (0, \delta']\;,
\ee
{\bf Fix $\delta'$ with this property}. Set
\be
\ep_1 = \min\left\{ \;\frac{1}{32 C_0 }\;,\; c_1\;,\;  \frac{1}{4E}\;,\;
\frac{1}{\hd\, \rho^{p_0+2} } \; , \; \frac{c_0r_0}{\gamma_1^{n_1}}\; , \; \frac{c_0^2(\gamma-1)}{16T\gamma_1^{n_1}}\, \right\} \;,
\ee
and let $b\in \R$ be such that $|b| \geq 1$ and
\begin{equation}
 \frac{\ep_1}{|b|} \leq \delta' \;.
\end{equation}

Let $\cc_m = \cc_m^{(\ep_1/|b|)}$ ($1\leq m \leq p$) be the family of {\it maximal  closed cylinders} 
contained in $\overline{U_0}$ with  $\diam(\cc_m )\leq \frac{\ep_1}{|b|}$ such that 
$U_0 \subset \cup_{j=m}^p \cc_m$ and $\overline{U_0} = \cup_{m=1}^p \cc_m$ (see Definitions 4.1).
It follows from (5.5), (5.4) and Lemma 3.2, that the length of each $\cc_m$ is not less than $n_1$, so
$\sigma^{n_1}$ is expanding on $\cc_m$.
Moreover, Proposition 3.3(a) implies that  $\diam(\cc_m)\geq \rho\, \frac{\ep_1}{|b|}$ for all $m$, so
\be
 \rho\, \frac{\ep_1}{|b|} \leq \diam(\cc_m ) \leq  \frac{\ep_1}{|b|} \quad, \quad 1\leq m\leq p\;.
\ee

Fix an integer $q_0 \geq 1$ such that
\be
\theta_0 < \theta_1  - 32\, \rho^{q_0-1}\;.
\ee
Next, let $\dd_1, \ldots, \dd_q$  be the list of all closed cylinders contained in $\overline{U_0}$ 
that are {\it subcylinders of co-length} $p_0\, q_0$ of some $\cc_m$ ($1\leq m\leq p$). That is, 
if $k_m$ is the length of $\cc_m$,  we consider the subcylinders of length $k_m  + p_0\, q_0$ of 
$\cc_m$, and we do this for any $m = 1, \ldots,p$.
Then  $\overline{U_0} = \cc_1 \cup \ldots \cup \cc_p =  \dd_1 \cup \ldots \cup \dd_q\; .$
Moreover, it follows from the properties of $\cc_m$ and Proposition 3.3 that
\be
\rho^{p_0\, q_0+1}\cdot \frac{\ep_1}{|b|} \leq \diam(\dd_j ) \leq  \rho^{q_0}\cdot 
\frac{\ep_1}{|b|} \quad, \quad 1\leq j\leq q\;.
\ee

Given $j = 1, \ldots,q $, $\ell = 1, \ldots, \ell_0$  and $i = 1,2$, set
$Z_j = \overline{\sigma^{n_1}(\hdd_j)}$, $\xijl = \overline{\vl_i(\hZ_j)}$,
$\hZ_j = Z_j \cap \hU$, $\hdd_j = \dd_j \cap \hU$, 
and $\hxijl = \xijl \cap \hU$. 
It then follows that $\dd_j =  \psi(Z_j)$, 
and $U = \cup_{j=1}^q Z_j$. Moreover, by Proposition 3.1, all $\xijl$ are cylinders.

$$\def\normalbaselines{\baselineskip20pt
      \lineskip3pt \lineskiplimit 3pt}
   \def\mapright#1{\smash{
       \mathop{-\!\!\!-\!\!\!-\!\!\!-\!\!\!-\!\!\!\longrightarrow}\limits^{#1}}}
    \def\mapleftt#1{\smash{
        \mathop{-\!-\!-\!-\!\!\!\longleftarrow}\limits_{#1}}}
  \def\mapdown#1{\Big\downarrow\rlap
        {$\vcenter{\hbox{$\scriptstyle#1$}}$}}
\matrix{U_0 & 
\stackrel{\di\mathop{-\!\!\!-\!\!\!-\!\!\!-\!\!\!-\!\!\!\longrightarrow}\limits^{\sigma^{n_1}}}{
\mathop{\longleftarrow \!\!\!-\!\!\!-\!\!\!-\!\!\!-\!\!\!-\!\!\!-}\limits_{\psi} } 
& \uu & 
\stackrel{\di\mathop{-\!\!\!-\!\!\!-\!\!\!-\!\!\!-\!\!\!\longrightarrow}\limits^{\vl_i}}{
\mathop{\longleftarrow \!\!\!-\!\!\!-\!\!\!-\!\!\!-\!\!\!-}\limits_{\sigma^{N}} }
& \vl_i(\uu) & \subset U_1\cr
     \bigcup &             & \bigcup &              & \bigcup  &  \cr
     \dd_j   &  \mapright{\sigma^{n_1}} & Z_j = \sigma^{n_1}(\dd_j) &  
     \mapright{\vl_i}  & \xijl  & \cr}$$

\bs

\noindent
{\bf Remark 5.5.} It follows from (5.2) that $\xijl \cap {X}_{i', j'}^{(\ell')} = \e$  whenever $(i,j,\ell) \neq  (i',j', \ell')$.

\bs

By Lemma 5.2(b), the {\it characteristic function} 
$\omega_{i,j}^{(\ell)} = \chi_{\hxijl} : \hU \longrightarrow [0,1]$ of $\hxijl$
belongs to $\clip_D(\hU)$ and $\Lip_D(\eijl) \leq 1/\diam(\xijl)$. Since
$\sigma^{N}(\hxijl) = \hZ_j = \sigma^{n_1}(\hdd_j)$ and $\sigma^{N}$ is expanding and one-to-one on $\hxijl$, 
it follows that $\sigma^{N - n_1}(\hxijl) = \hdd_j$ and (2.1) gives
$\diam(\dd_j) \leq \frac{\gamma_1^{N-n_1}}{c_0}\, \diam(\xijl) 
\leq \frac{\gamma_1^{N}}{c_0}\, \diam(\xijl)\;.$
Combining this with (5.8) gives
\be
\diam(\xijl) \geq \frac{c_0\, \rho^{p_0\, q_0+1}}{\gamma_1^{N}} \cdot \frac{\ep_1}{|b|}
\ee
for all $i = 1,2$, $j = 1, \ldots,q $ and $\ell = 1, \ldots, \ell_0$.

Let $J$  be a {\it subset of the set} 
$\Xi = \Xi(a,b) = \{\; (i,j, \ell) \; :  \; 1\leq i \leq 2\; ,\;  1\leq j\leq q\; , \;  1\leq \ell \leq \ell_0\;\}.$
Set
\be
\mu = \mu(N) = \min \left\{\; \frac{1}{4} \; , \;  
\frac{c_0  \, \rho^{p_0q_0+2}\, \ep_1}{4\, \gamma_1^N }\; , \;
\frac{1}{4\,e^{2 T N}} \,\sin^2\left(\frac{\hd\, \rho\, \ep_1}{256}\right) \; \right\}\;,
\ee
and define the function  
$\beta = \beta_{J} : \hU \longrightarrow [0,1]$ by
$\di \beta  = 1- \mu \,\sum_{(i, j, \ell) \in J} \eijl\;.$
Clearly $\beta \in \clip_D(\hU)$ and $1-\mu \leq \beta(u) \leq 1$ for any $u \in \hU$.
Using Remark 5.5, Lemma 5.2(b)  and (5.9) one derives that
\begin{equation}
\Lip_{D} (\beta) \leq \Gamma =  \frac{ 2 \mu \,\gamma_1^N}{c_0\,\rho^{p_0q_0+2}}\cdot \frac{|b|}{\ep_1} \; .
\end{equation}

Next, define {\it the operator} $\nn = \nn_J(a,b) : \clip_D (\hU) \longrightarrow \clip_D (\hU)$ by 
$\left(\nn h\right) = \ma^N (\beta_J \cdot h)$.

\ms

The following lemma contains statements similar to 
Proposition 6 and Lemma 11 in \cite{kn:D2} and by means of Lemma 5.4 their proofs are
also very similar, so we omit them.

\bs

\noindent
{\bf Lemma 5.6.} 
{\it Under the above conditions for $N$ and $\mu$ the following hold :}

(a) $\nn h\in K_{E|b|}(\hU)$ {\it for any} $h\in K_{E|b|}(\hU)$;

(b) {\it If $h \in \clip_D (\hU)$ and $H \in K_{E|b|}(\hU)$ are such that $|h|\leq H$ in $\hU$ and
$| h (v) - h(v')|\leq E|b| H(v')\, D (v,v')$
for any $v,v'\in U_j$, $j = 1,\ldots,k$, then for any $i = 1, \ldots,k$ and any $u,u'\in \hU_i $ we have
$| (\lab^N  h)(u) - (\lab^N h)(u')| \leq E |b| (\nn H)(u')\, D (u,u')$.}

\bs

\noindent
{\bf Definition.} Given $t > 0$ and $S > 0$, a subset $W$ of $\hU$ will be called 
{\it $(t,S)$-dense} in $\hU$ if for every $u\in \hU $ 
there exist a cylinder $\cc$ containing $u$ with $\diam(\cc) \leq S t$ and a cylinder $\cc'$ 
with $\diam(\cc') \geq t$ such that  $\hc' \subset W\cap \cc$.

\ms
Below we use  the constants $c_1,c_2$ from (3.1) and $\|g\|_0 = \sup_{u\in \hU} |g(u)|$.

\bs

\noindent
{\bf Lemma 5.7.} 
{\it Let $A > 0$, $S \geq 1$ and let $\di \epsilon = \epsilon(S,A) = \frac{d_S}{ e^{S^2 A^2} }$,
where $d_S = \frac{c_1}{c_2}\, e^{-p_0\, \|g\|_0\, \left( \frac{\ln S}{|\ln \rho |} + 1\right)}$.
Then for any $t > 0$, any $(t,S)$-dense subset
$W$ of $\hU$ and any $H \in K_{A/t}(\hU)$, we have}
$$\int_W H^2 \; d\nu \geq \epsilon \int_{\hU} H^2 \; d\nu\; .$$

\ms

\noindent
{\it Proof of Lemma 5.7.}  Let $t > 0$, let $W$ be a $(t,S)$-dense subset of $\hU $ and let 
$H \in  K_{A/t}(\hU)$. Let $\bb_1, \ldots, \bb_n$ be the maximal cylinders in $U$ with 
$\diam(\bb_j) \leq S t$ for any $j = 1, \ldots, n$.
Then $\cup_{j=1}^n \hbb_j = \hU$, and so $\sum_{j=1}^n \nu(\hbb_j) = 1$.
Setting $m_j = \inf_{u\in \bb_j } H(u)$ and  $M_j = \sup_{u\in \bb_j } H(u)$,
it  follows from $H \in K_{A/t}(\hU)$ that $\frac{M_j}{m_j} \leq e^{\frac{A S t}{t}} = e^{A S}$.
Thus, $H(u) \geq m_j \geq M_j e^{-  A S}$ for all $u\in \hbb_j$.

For each $j = 1, \ldots,n$ choose an arbitrary point  $u_j\in \hbb_j$. Since $W$ is $(t,S)$-dense
in $\hU$, there exists a subcylinder $\bb'_j$ of $\bb_j$ such that $\diam(\bb'_j) \geq t$
and $\hbb'_j \subset W\cap \bb_j$. If $q_j$ is the co-length of $\bb'_j$ in $\bb_j$, and
$q_j = r_jp_0 + s_j$ for some integers $r_j \geq 0$, $0\leq s_j < p_0$, by Proposition 3.3(b) we have
$t \leq \diam(\bb'_j) \leq \rho^{r_j}\, \diam(\bb_j) \leq \rho^{r_j}\, S t\;,$
so $\rho^{r_j} \geq 1/S$, i.e. $r_j \leq \ln S/ |\ln \rho|$. Thus, 
$q_j < p_0(r_j+1) \leq p_0\, (\ln S/|\ln \rho| + 1)$. If $p_j$ is the length of
the cylinder $\bb_j$ and $p'_j = p_j+q_j$ that of $\bb'_j$,   (3.1) gives
$\frac{\nu(\bb'_j)}{\nu(\bb_j)} \geq \frac{c_1}{c_2}\, e^{-q_j\, \|g\|_0} \geq d_S\;.$
Hence $\nu(\bb'_j) \geq d_S\, \nu(\bb_j)$ for all $j = 1,\ldots,n$.

It now follows that $\nu(W\cap \bb_j) \geq \nu(\bb'_j) \geq d_S\, \nu(\bb_j)$ for any $j$, so
\begin{eqnarray*}
\int_W H^2(u) \; d\nu(u)
& =     & \sum_{j=1}^n \int_{W \cap \bb_j} H^2(u) \; d\nu(u)
 \geq   \sum_{j=1}^n M^2_j e^{-A^2 S^2}  \nu(W \cap \bb_j)\\
& \geq & e^{-A^2 S^2} \sum_{j=1}^n M^2_j\, d_{S} \nu(\bb_j)
 \geq    \frac{d_S}{ e^{A^2 S^2}}\,  \int_{\hU} H^2(u) \; d\nu(u)
= \epsilon  \int_{\hU} H^2(u) \; d\nu(u) \; .
\end{eqnarray*}
This proves the assertion. \endofproof\\

\noindent
{\bf Definitions.} A subset $J$ of $\Xi$ will be called {\it dense} if for any $m = 1,\ldots, p$ 
there exists $(i,j, \ell)\in J$ such that $\dd_j \subset \cc_m$.
Denote by $\J = \J(a,b)$ the {\it set of all dense subsets} $J$ of $\Xi$. 

\bs

\noindent
{\bf Lemma 5.8.} 
{\it Given the number $N$, there exist  $\epsilon_2 = \epsilon_2(N) \in (0,1)$ and $a_0 = a_0(N) > 0$ 
such that $\di\int_{\hU} (\nn H)^2d\nu \leq (1-\epsilon_{2}) \,\int_{\hU} H^2 d\nu\; $ whenever 
$|a| \leq a_0$, $J$ is dense and $H \in K_{E|b|}(\hU)$.
More precisely, we can take $\epsilon_2 = \frac{\epsilon'\,\mu \,e^{-NT}}{4}$ and  
$a_0  =  \min \left\{ a'_0 , \frac{1}{C_0\, N}\, \ln\left(1+  \frac{\epsilon'\,\mu 
\,e^{-NT}}{4}\right)  \right\}\;, $
where $ \epsilon' = \frac{d_{S}}{e^{ (E \, \gamma_1^{n_1}\ep_1/c_0)^2}}$ and
$S = \frac{\gamma_1^{n_1}}{c^2_0\, \gamma^{n_1} \rho^{p_0\, q_0 + 1}}$. }

\bs

\noindent
{\it Proof of Lemma} 5.8. The definition of $\nn$ and the Cauchy-Schwartz  inequality imply
\be
(\nn H)^2 (u) = (\ma^{N} (\beta H))^2(u)\leq (\ma^{N} \beta^2)(u)\cdot (\ma^{N} H^2)(u)
\leq (\ma^{N} H^2)(u)
\ee
for all $u \in \hU$.

Denote $W =  \cup_{(i,j,\ell)\in J} \hZ_j$. Then  $u\in W$ means that there exists $(i,j, \ell) \in J$
with $\vl_i(u) \in \xijl$, and so  $\beta(\vl_i(u)) = 1 -\mu$.

We will now show that $W$ is $(t,S)$-dense in $\hU$, where
$t = c_0\, \gamma^{n_1}\, \rho^{p_0\, q_0 + 1 } \cdot \frac{\ep_1}{|b|}$ and $S$  is as in the statement of the lemma.
Let $u\in \hU$. Since $\hU \subset \uu \subset \cup_{m=1}^p \sigma^{n_1}(\cc_m)$, we have
$u \in \cc = \sigma^{n_1}(\cc_m)$ for some $m$. 
Since $J$ is dense, there exists $(i,j,\ell) \in J$ so that $\dd_j \subset \cc_m$.
Then $Z_j = \sigma^{n_1}(\dd_{j}) \subset \cc$, so $\hZ_j \subset W \cap \cc$.
Now  (2.1), (5.6) and (5.8) yield $\diam(\hZ_j) \geq t$ and $\diam(\cc) \leq St$, so
$W$ is $(t,S)$-dense in $\hU$.

Let $H \in K_{E|b|}(\hU)$.  Setting $A = E  \,c_0\, \gamma^{n_1}\, \rho^{p_0\, q_0 + 1}\, \epsilon_1$, 
we have $H \in  K_{\frac{A}{t}}(\hU )$. Since $\epsilon' = \frac{d_{S}}{e^{(S A)^2}}$, it follows that 
$\ep' = \epsilon(S,A)$, the number defined in Lemma 5.7. By Lemma 5.6 (a), $\nn H\in K_{A/t}(\hU)$, so 
Lemma 5.7 implies
\begin{equation}
\int_W (\nn H)^2 \; d\nu \geq \epsilon' \int_{\hU} (\nn H)^2 \; d\nu \;.
\end{equation}

It follows from the definition of $\beta$ that $\ma^{N}\beta^2 \leq  1 - \mu \, e^{-NT}$ on $W$.
Using this,  (5.12) and  (5.13), as in \cite{kn:D2} (see also the proof of Lemma 7.4 in \cite{kn:St1}) we get
$$\int_{\hU} (\nn H)^2 \, d\nu \leq \int_{\hU} (\ma^N H^2)\, d\nu 
- \mu \, e^{-NT} \,\ep'\,  \int_{\hU} (\nn H)^2\, d\nu\;,$$
and therefore
$\di \int_{\hU} (\nn H^2)\; d\nu \leq   \frac{1}{1+ \mu \, \ep'\, e^{-NT}}  
\int_{\hU} (\ma^N H^2)\; d\nu \;.$
Since 
$\ma^N H^2 = L_{f^{(0)}}^N\left(e^{(\fa - f^{(0)})_N}\, H^2\right) \leq e^{N|a|C_0}\, L_{f^{(0)}}^N H^2$
on $\hU$, it now follows from $|a| \leq a_0$ and $L_{f^{(0)}}^* \nu = \nu$ (see section 2) that
$$\int_{\hU} (\nn H^2)\; d\nu \leq   \frac{e^{N|a|C_0}}{1+ \mu \, \ep'\, e^{-NT}}  
\int_{\hU} (L_{f^{(0)}}^N H^2)\; d\nu 
\leq   \frac{1+ \ep_2}{1+ 4\ep_2}  \int_{\hU} H^2\; d\nu \leq   (1-\ep_2)\,  \int_{\hU} H^2\; d\nu\;.$$
This completes the proof of the lemma.
\endofproof

\bs

In what follows we assume that $h, H\in \clip_D (\hU)$ are such that 
\begin{equation}
H\in K_{E|b|}(\hU)\quad , \quad |h(u)|\leq H(u) \:\:\:\:, \:\:\; u\in \hU \;,
\ee
and
\be
|h(u) - h(u')|\leq E|b| H(u')\, D (u,u')\:\:\:\:\:
\mbox{\rm whenever} \:\: u,u'\in \hU_i\;, \; i = 1, \ldots,k \;.
\end{equation}  

Define the functions $\geil : \hU  \longrightarrow \C$ ($\ell = 1, \ldots, j_0$, $i = 1,2$) by
$$\di \geol(u) = \frac{\di \left| e^{(\fa_{N}+\i b\tau_{N})(\vl_1(u))} h(\vl_1(u)) +
 e^{(\fa_{N}+\i b\tau_{N})(\vl_2(u))} h(\vl_2(u))\right|}{\di (1-\mu)
e^{\fa_{N}(\vl_1(u))}H(\vl_1(u)) + e^{\fa_{N}(\vl_2(u))}H(\vl_2(u))}\; ,$$ 
$$\di \getl(u) = \frac{\di \left| e^{(\fa_{N}+\i b\tau_{N})(\vl_1(u))} h(\vl_1(u)) +
 e^{(\fa_{N}+\i b\tau_{N})(\vl_2(u))} h(\vl_2(u))\right|}{\di
e^{\fa_{N}(\vl_1(u))}H(\vl_1(u)) +  (1-\mu) e^{\fa_{N}(\vl_2(u))}H(\vl_2(u))}\; ,$$ 
and set $ \gl(u) = b\, [\tau_N(\vl_2(u)) - \tau_N(\vl_1(u))]$, $u\in \hU \;.$
 
\bs

\noindent
{\bf Remark.} It is easy to see that for any $j$ and $\ell$ the set
$\{ \gl(u) : u \in Z_j\}$ is contained in an interval of length $< 1/8$. 
Indeed, given $i$, $j$ and $\ell$ and $u, u'\in \hZ_j$, consider
$x = \vl_i(u), x' = \vl_i(u') \in \hxijl$. Since $\sigma^{N-n_1}(\xijl) = \dd_j$, by (5.8) and (2.1),
$\diam(\sigma^N(\xijl)) = \diam (\sigma^{n_1}(\dd_j)) \leq \frac{\gamma_1^{n_1}\ep_1}{c_0\, |b|}$,
so $d(u,u') = d(\sigma^N(x), \sigma^N(x')) \leq \frac{\gamma_1^{n_1}\ep_1}{c_0\, |b|}$. This and (2.1)
give $d(\sigma^j(x), \sigma^j(x')) \leq \frac{1}{c_0\, \gamma^{N-j}}\, d(u,u') 
\leq \frac{\gamma_1^{n_1}\ep_1}{c^2_0 \gamma^{N-j} |b|}$, so by (5.4),
$|\tau_N(v) - \tau_N(v')| \leq \sum_{j=0}^{N-1}|\tau(\sigma^j(x))-\tau(\sigma^j(x'))|
\leq  \frac{T\, \gamma_1^{n_1}\,\ep_1}{c^2_0\, |b|\, (\gamma-1)} < \frac{1}{16\, |b|}$. 
Thus, the set $\{ \tau_N(\vl_i(u)) : u\in Z_j\}$
is contained in an interval of length $< \frac{1}{16|b|}$, and therefore\\
$\{ \gl(u) : u \in Z_j\}$ is contained in an interval of length $< 1/8$. 

\bs

\noindent
{\bf Definitions.} 
We will say that the cylinders $\dd_j$ and $\dd_{j'}$ are {\it adjacent} if they 
are subcylinders of the same $\cc_m$ for some $m$. 
If $\dd_j$ and $\dd_{j'}$ are contained in $\cc_m$  for some $m$
and  for some  $\ell = 1, \ldots, \ell_0$ there exist  $u \in \dd_j$ and $v\in \dd_{j'}$  
such that $d(u,v) \geq \frac{1}{2}\, \diam(\cc_m)$ and 
$\left\la \frac{r^{-1}(v) - r^{-1}(u)}{\| r^{-1}(v) - r^{-1}(u)\|}\;,\; \eta_\ell \right\ra  \geq \theta_1$,
we will say that $\dd_j$ and $\dd_{j'}$ are {\it $\eta_\ell$-separable in $\cc_m$}.

\bs

\noindent
{\bf Lemma 5.9.} {\it Let $j, j'\in \{ 1, 2,\ldots,q\}$ be
such that $\dd_j$ and $\dd_{j'}$ are contained in $\cc_m$ and  are $\eta_\ell$-separable 
in $\cc_m$ for some $m = 1, \ldots, p$ and $\ell = 1, \ldots, \ell_0$ . Then
$ |\gl (u) - \gl(u')| \geq c_2\epsilon_{1}$ for all $u\in \hZ_j$ and  $u'\in \hZ_{j'}$,
where $\di c_2 =  \frac{\hd\, \rho}{16}$.} 

\bs

\noindent
{\it Proof of Lemma 5.9.} Let $u\in \hZ_j$ and $u'\in \hZ_{j'}$; then  $x = \psi(u) \in \hdd_j$ and
$x' = \psi(u') \in \hdd_{j'}$. Also  $x = r(s)$ and $x' = r(s')$ for some $s, s'\in V_0$. Set 
$\eta = \frac{s- s'}{\|s-s'\|} \in \S^{n-1}\;.$

Since $\dd_j$ and $\dd_{j'}$ are $\eta_\ell$-separable in $\cc_m$, there exist
$x_0 = r(s_0) \in \dd_j$ and $x'_0 = r(s'_0) \in \dd_{j'}$  such that 
$d(x_0,x'_0) \geq \frac{1}{2}\, \diam(\cc_m)$ and
$\left\la \eta_0 \;,\; \eta_\ell \right\ra  \geq \theta_1$,
where $\di \eta_0 = \frac{s_0 - s'_0}{\|s_0 - s'_0\|} \in \S^{n-1}\;.$ By (4.2), (5.8) and (5.6),
$\|s-s_0\| \leq 2\, d(r(s), r(s_0)) \leq 2\, \diam(\dd_j) \leq 2\rho^{q_0-1}\, \diam(\cc_m)\;,$
and similarly $\|s'-s'_0\| \leq 2\rho^{q_0-1}\, \diam(\cc_m)$. This implies
\be
\left| \|s-s'\| - \|s_0-s'_0\|\right| \leq \|s-s_0\| + \|s'-s'_0\| \leq 4\rho^{q_0-1} \diam(\cc_m)\;.
\ee
Hence
$\|\eta_0 - \eta\| =   \left\| \frac{s_0 - s'_0}{\| s_0 - s'_0\|} - \frac{s - s'}{\| s - s'\|} 
\right\| \leq  \frac{ 8\rho^{q_0-1} \diam(\cc_m)}{\|s_0-s'_0\|} 
 \leq  \frac{ 16 \rho^{q_0-1} \diam(\cc_m)}{d(x_0,x'_0)} \leq 32 \rho^{q_0-1}$.
Combining this with  (5.7) gives
$ \la \eta , \eta_\ell\ra =  \la \eta_0 , \eta_\ell\ra + \la \eta - \eta_0 , \eta_\ell\ra
\geq \theta_1 - 32 \rho^{q_0-1} > \theta_0\;.$
Thus, $\eta \in B_\ell$, and Lemma 4.2 implies
$\left| \left[I_{\eta , h} \left(\tau_{N}(\vl_2(\trr(\cdot))) 
- \tau_{N}(\vl_1(\trr(\cdot)))\right)\right](\hs) \right| \geq \frac{\hd}{2} $
for all $\hs\in r^{-1}(U_0)$ and $h\neq 0$ such that 
$\hs+ h \eta \in r^{-1}(U_0\cap \mt)$. 

Since $u = \sigma^{n_1}x = \tr(s)$ and $u' = \tr(s')$, we have $s, s'\in r^{-1}(U_0\cap \mt)$ and 
$s = s' + h \eta$ 
with $h = \| s - s'\|$.  It then follows from the above, (5.16) and (5.6) that
\begin{eqnarray*}
\frac{1}{|b|}\, |\gl(u) - \gl(u') | 
&      =  & \left| [\tau_{N}(\vl_2(\trr(s))) - \tau_{N}(\vl_1(\trr(s)))] 
-  [\tau_{N}(\vl_2(\trr(s'))) - \tau_{N}(\vl_1(\trr(s')))]  \right|\\
&\geq  & \frac{\hd}{2} \, \|s-  s'\| \geq 
 \frac{\hd}{2} \left( \frac{1}{4}  \diam(\cc_m) - 4\rho^{q_0-1}\diam(\cc_m) \right)  \geq 
\frac{\hd\,  \rho\, \ep_1}{16 \,|b|}\;. \:\: \endofproof
\end{eqnarray*}

\ms

The  following lemma is central for this section.

\bs

\noindent
{\bf Lemma 5.10.}  {\it Assume $b$ is chosen in such a way that {\rm (5.5)} holds. Then  for any
$j = 1, \ldots,q$  there exist $i \in \{ 1,2\}$, $j' \in \{ 1,\ldots,q\}$ and 
$\ell \in \{ 1, \ldots, \ell_0\}$  such that  $\dd_j$ and $\dd_{j'}$ 
are adjacent and $\chi^{(i)}_{\ell} (u) \leq 1$ for all  $u\in \hZ_{j'}$} .

\bigskip

To prove this  we need the following lemma which
coincides with Lemma 14 in \cite{kn:D2} and its proof is almost the same, so we omit it.

\bs

\noindent
{\bf Lemma 5.11.}  {\it If $h$ and $H$ satisfy {\rm (5.14)-(5.15)},
then for any $j = 1, \ldots,q$, $i = 1,2$ and $\ell = 1,\ldots,  \ell_0$ we have:}

(a) {\it $\di\frac{1}{2} \leq \frac{H(\vl_i(u'))}{H(\vl_i(u''))} \leq 2$ for all}
$u', u'' \in \hZ_j$;

(b) {\it Either for all $u\in \hZ_j$ we have
$|h(\vl_i(u))|\leq \frac{3}{4}H(\vl_i(u))$, or  $|h(\vl_i(u))|\geq \frac{1}{4}H(\vl_i(u))$
for all $u\in \hZ_j$.} 

\bs

\noindent
{\it Proof of Lemma 5.10.} Given $j = 1, \ldots, q$, let $m = 1, \ldots, p$
be such that $\dd_j \subset \cc_m$.  By (5.5), $\delta = \ep_1/|b| \in (0,\delta']$, so
it follows from (5.3) that $\cc_m \subset M^{(\delta)}_{\eta_\ell}(U_0)$
for some $\ell = 1, \ldots, \ell_0$. This means that there exist $u,v\in \cc_m$ such that 
$d(u,v) \geq \frac{1}{2}\, \diam(\cc_m)$ and
$\left\langle \frac{r^{-1}(v) - r^{-1}(u)}{\|r^{-1}(v) - r^{-1}(u)\|}\; , \; \eta_\ell\; \right\rangle 
 \geq \theta_1\;.$
Let $j',j'' = 1, \ldots,q$ be such that $u\in \dd_{j'}$ and $v\in \dd_{j''}$. 
(Notice that we may have $j' = j$ or $j'' = j$.) Then
$\dd_{j'}$ and $\dd_{j''}$ are $\eta_\ell$-separable in $\cc_m$.

Fix $\ell$, $j'$ and $j''$ with the above properties, and set
$\hZ = \hZ_j \cup \hZ_{j'}\cup \hZ_{j''}\; .$ 
If there exist $t \in \{j, j', j''\}$ and $i = 1,2$ such 
that the first alternative in Lemma 5.11(b) holds for $\hZ_{t}$, $\ell$  and $i$, then $\mu \leq 1/4$
implies $\chi_\ell^{(i)}(u) \leq 1$ for any $u\in \hZ_{t}$.

Assume that for every $t\in \{ j, j', j''\}$ and every $i = 1,2$ the second alternative 
in Lemma 5.11(b) holds for $\hZ_{t}$, $\ell$ and $i$, i.e. 
$|h(\vl_i(u))|\geq \frac{1}{4}\, H(\vl_i(u))$, $u \in \hZ$.

Since $\psi(\hZ) = \hdd_j \cup \hdd_{j'} \cup \hdd_{j''} \subset \cc_m$, given $u,u'\in \hZ$ we have 
$\sigma^{N-n_1}(\vl_i(u)), \sigma^{N-n_1}(\vl_i(u'))\in \cc_m$. 
Notice that by Proposition 3.1, $\cc' = \vl_i\circ \wl_i(\cc_m)$ is a cylinder. Moreover,
$\sigma^{N-n_1}$ is expanding on $\cc'$, so by (2.1) and (5.6),
$\diam(\cc') \leq \frac{\diam(\cc_m)}{c_0\, \gamma^{N-n_1}} \leq 
\frac{\ep_1}{c_0\, \gamma^{N-n_1}\, |b|}\;.$
Thus, using the above assumption, (5.14),  (5.15),  (5.1) and the constant $c_2$ from Lemma 5.9, 
and assuming e.g.   $|h(\vl_i(u))| \geq |h(\vl_i(u'))|$, we get
\begin{eqnarray*}
\frac{|h(\vl_i(u)) - h(\vl_i(u'))|}{\min\{ |h(\vl_i(u))| , |h(\vl_i(u'))| \}}
 \leq  \frac{E|b|\, H(\vl_i(u'))}{|h(\vl_i(u'))| } D (\vl_i(u),\vl_i(u'))
 \leq  4 E|b|\,  \diam(\cc')
 \leq \frac{c_2 \ep_1}{8} < \frac{1}{2}\;.
\end{eqnarray*}
Thus, the angle between the complex numbers  $h(\vl_i(u))$ and $ h(\vl_i(u'))$ 
(regarded as vectors in $\R^2$)  is less than $\pi/3$. 
In particular,  for any $i = 1,2$ we can choose real continuous functions $\theta_i(u)$, $u \in  \hZ$,
with values in $[0,\pi/3]$ and constants $\lambda_i \in [0,2\pi)$  such that
$\di h(\vl_i(u)) = e^{\i(\lambda_i +\theta_i(u))}|h(\vl_i(u))|$ for all $u\in \hZ$.
Using the above, $\theta \leq 2 \sin \theta$ for $\theta \in [0,\pi/3]$, and some elementary geometry  yields
$|\theta_i(u) - \theta_i(u')|\leq 2 \sin |\theta_i(u) - \theta_i(u')| < \frac{c_2\ep_1}{8}\;.$

The difference between the arguments of the complex numbers
$e^{\i \,b\,\tau_N(\vl_1(u))} h(\vl_1(u))$ and $e^{\i \,b\, \tau_N(\vl_2(u))} h(\vl_2(u))$
is given by the function
$$\Gl(u) = [b\,\tau_N(\vl_1(u)) + \theta_1(u) + \lambda_1] -  [b\, \tau_N(\vl_2(u)) + \theta_2(u) + \lambda_2]\;.$$
Notice that by the Remark before Lemma 5.9, the set 
$\{ |\gl(u') - \gl(u'')| : u'\in \hZ_{j'}, u''\in \hZ_{j''}\}$ is contained in an
interval of length $< 1/4$.
Given $u'\in \hZ_{j'}$ and $u''\in \hZ_{j''}$, since $\hdd_{j'}$ and $\hdd_{j''}$ are 
contained in $\cc_m$ and are $\eta_\ell$-separable in $\cc_m$, it follows from  Lemma 5.9 and the above that
\begin{eqnarray*}
|\Gl(u')- \Gl(u'')| \geq  |\gl(u') - \gl(u'')| - |\theta_1(u')-\theta_1(u'')| 
- |\theta_2(u')-\theta_2(u'')| \geq   \frac{c_2\ep_1}{2}\;.
\end{eqnarray*}
Thus,  $|\Gl(u')- \Gl(u'')|\geq \frac{c_2}{2} \epsilon_{1}$ for all $u'\in \hZ_{j'}$ and 
$u''\in \hZ_{j''}$. Hence either 
$|\Gl(u')| \geq \frac{c_2}{4}\epsilon_{1}$ for all $u'\in \hZ_{j'}$ or 
$|\Gl(u'')| \geq \frac{c_2}{4}\epsilon_{1}$ for all $u''\in \hZ_{j''}$.
Assuming for example that $|\Gl(u)| \geq \frac{c_2}{4}\epsilon_{1}$ for all 
$u\in \hZ_{j'}$, as in \cite{kn:D2} (see also \cite{kn:St1}) one shows that
$\chi_\ell^{(1)} (u)  \leq 1$ and $\chi_\ell^{(2)} (u)  \leq 1$ for all $u \in \hZ_{j'}$.  \endofproof

\bs

\noindent
{\it Proof of Lemma 5.3.} As before, we assume that $N$  and $\mu > 0$ 
satisfy (5.1) and (5.10). Define $\ep_1$ by (5.4), take $\ep_2 = \ep_2(N) > 0$ and $a_0 = a_0(N) > 0$
as in Lemma 5.8 and set $\hat{\rho} = 1- \ep_2 $.

Let $a\in \R$ and $b \in \R$ be such that $|a| \leq a_0$ and $|b| \geq b_0 = \ep_1/\delta'$.
Then for any $J \in \J(a,b)$, Lemma 5.6 (a) implies property (a) in Lemma 5.3 for the
operator $\nn_J$,  while Lemma 5.8 gives property (b) in Lemma 5.3.

To check (c) in Lemma 5.3, assume that $h, H \in  \clip_D (\hU)$ satisfy (5.14) and (5.15). 
Now define the subset $J$ of $\J(a,b)$ in the following way. First, include in $J$ all $(1,j, \ell)\in \Xi$ such
that $\geol(u) \leq 1$ for all $u \in \hZ_j$. Then for any $j = 1, \ldots, q$ and $\ell = 1, \ldots, j_0$ include
$(2,j,\ell)$ in $J$ if and only if $(1,j,\ell)$ has not been included in $J$ (that is, $\geol(u) > 1$ for 
some $u \in \hZ_j$) and $\getl(u) \leq 1$ for all $u \in \hZ_j$. It follows from Lemma 5.10   that $J$ is dense.  
(Clearly, $J$ depends not only on $N$, $a$ and $b$, but on $h$ and $H$ as well.)

Consider the operator $\nn = \nn_{J}(a,b) : \clip_D (\hU) \longrightarrow \clip_D (\hU)$. 
Then Lemma 5.6 (b) implies $| (\lab^N  h)(u) - (\lab^N h)(u')| \leq E |b| (\nn H)(u')\, D (u,u')$
whenever $u,u'\in \hU_i $ for some $i = 1, \ldots,k$.  So, it remains to show that
\begin{equation}
\left| (\lab^N h)(u)\right| \leq (\nn H)(u) \:\:\: , \:\: u\in \hU \;.\
\end{equation}

Let $u \in \hU$.  If $u \notin \hZ_j$ for any $(i,j, \ell)\in J$, then $\beta (v) = 1$ whenever
$\sigma^N v = u$ and therefore $\left| (\lab^N h)(u)\right|  \leq  (\ma^N (\beta H))(u) = (\nn H)(u)$.

Assume that $u\in \hZ_j$  e.g. for  $(1,j, \ell)\in J$; then $(2,j, \ell) \notin J$. Since $\geol(u) \leq 1$, 
$\beta (\vl_1(u)) \geq 1-\mu$ and $\beta (\vl_2(u)) = 1$, using (5.14) one derives
\begin{eqnarray*}
\left| (\lab^N h)(u)\right| 
&\leq  & \sum_{\sigma^N v = u, \;v\neq v_1(u),v_2(u)} e^{\fa_N(v)} |h(v)|\\
&      &  + \left[e^{\fa_N(v_1(u))} \beta (v_1(u)) H(v_1(u)) + 
 e^{\fa_N(v_2(u))} \beta (v_2(u))H(v_2(u))\right] \leq (\nn H)(u)\;,
\end{eqnarray*}
which proves (5.17).
This completes the proof of Lemma 5.3. \endofproof

\def\Lye{L_{y,\eta}}
\def\Lyep{L^{(p)}_{y,\eta}}
\def\Fyp{F^{(p)}_y}
\def\Fxp{F^{(p)}_x}
\def\Lxx{L_{x,\xi}}
\def\Lxxp{L^{(p)}_{x,\xi}}
\def\chBo{\check{B}^{u,1}}
\def\tBo{\tB^{u,1}}
\def\hBo{\hB^{u,1}}
\def\hpi{\hat{\pi}}
\def\Bmt{\overline{B_{\ep_0}(\mt)}}

\section{Non-integrability conditions}
\setcounter{equation}{0}

Throughout we assume that $\phi_t$ is a $C^2$ contact flow on $M$ with a $C^2$ invariant  contact
form $\omega$. The following condition says  that $d\omega$ is in some sense 
non-degenerate on the `tangent space' of $\mt$ near some of its points:

 \ms

\noindent
{\sc (ND)}:  {\it There exist $z_0\in \mt$, $\ep_0 > 0$  and  $\mu_0 > 0$ such that
for any $\ep \in (0,\ep_0]$, any $\hz \in \mt \cap W^u_{\ep}(z_0)$ and any unit vector 
$\eta \in E^u(\hz)$ tangent to 
$\mt$ at $\hz$ there exist $\tz \in \mt \cap W^u_{\ep}(\hz)$, $\ty \in W^s_\ep(\tz)$ and 
a unit vector $\xi \in E^s(\ty)$ tangent to $\mt$ at $\ty$ with 
$|d\omega_{\tz}(\xi_{\tz},\eta_{\tz}) | \geq \mu_0$,
where $\eta_{\tz}$ is the parallel translate of $\eta$ along the geodesic in $W^u_{\ep}(\tz)$ from
$\hz$ to $\tz$, while $\xi_{\tz}$ is the parallel translate of $\xi$ along the geodesic in 
$W^s_{\ep}(\tz)$ from $\ty$ to $\tz$.}

\ms

\noindent
{\bf Remark.} It appears the above condition would become significantly more
restrictive if one requires the existence of a unit vector $\xi \in E^s(\tz)$ tangent to $\mt$ at 
$\tz$ with  $|d\omega_{\tz}(\xi,\eta_{\tz}) | \geq \mu_0$. The reason for this is that in general the set
of unit tangent vectors to $\mt$ does not have to be closed in the bundle $E^s_\mt$.
That is, there may exist a point  $\tz \in \mt$, a
sequence $\{ z_m\} \subset W^s_{\ep}(\tz)\cap \mt$ and for each $m$ a unit vector
$\xi_m$ tangent to $\mt$ at $z_m$ such that $z_m \to z$ and $\xi_m \to \xi$ as $m \to \infty$,
however $\xi$ is not tangent to $\mt$ at $\tz$. 

\ms

\noindent
{\bf Proposition 6.1.} {\it For contact flows $\phi_t$ with Lipschitz local stable holonomy maps, the condition {\rm (ND)} implies {\rm (LNIC)}.}

\ms

Since  (ND) is always satisfied when $\dim(M) =3$ or $\mt = M$, we get the following.

\ms

\noindent
{\bf Corollary 6.2.} {\it For contact flows $\phi_t$ with either $\dim(M) = 3$ (and an arbitrary basic set $\mt$)
or $\mt = M$ the condition {\rm (LNIC)}  is always satisfied on $\mt$.} \endofproof

\ms

To prove Lemma 6.1 we will make use of the following lemma which is a consequence of
Lemma B.7 in \cite{kn:L2}.

\ms

\def\tde{\tilde{\delta}}

\noindent
{\bf Lemma 6.3.} (\cite{kn:L2}) {\it Let $\phi_t$ be a contact flow on $M$ with a contact form $\omega$
and let $\mt$ be  a basic set for $\phi_t$ with Lipschitz local (un)stable holonomy maps. Then 
for every $\tz \in \mt$ and every $\tde > 0$ there exists $\tep \in (0,\ep_0)$ such that
for any $z\in \mt \cap W^u_{\tep}(\tz)$, any $x\in W^u_{\tep}(z) \cap \mt$ and any $y\in W^s_{\tep}(z) \cap \mt$ we have
$|\Delta(x,y)  - d\omega_z(u,v)| \leq \tde\, \|u\|\, \|v\|\;,$
where $u \in E^u(z)$ and $v \in E^s(z)$ are such that $\exp^u_z(u) = x$ and $\exp^s_z(v) = y$.}

\bs

Fix an arbitrary constant $\theta_0$  such that
\be
1- \frac{\mu^2_0}{128 C^2} \leq \theta_0  < 1 \;.
\ee

\noindent
{\it Proof of Proposition 6.1.} Assume that $\mt$ is a basic set for a contact flow $\phi_t$ such that 
(ND) holds on $\mt$.
Let $z_0\in \mt$, $1> \ep_0  > 0$ and $\mu_0 > 0$  be as in the statement of (ND) above. Fix a 
constant $C > 0$ with
$|d\omega_z(u,v)| \leq C\, \|u\|\, \|v\|$ for all $z\in \mt$ and $u,v \in T_zM$.

To check (LNIC), consider arbitrary $\ep \in (0,\ep_0]$ and $\hz\in \mt\cap W^u_\ep(z_0)$, and let 
$\eta \in E^u(\hz)$ be 
an arbitrary tangent vector to $\mt$ at $\hz$ with $\|\eta\| = 1$.
It follows from the condition (ND) that there exist $\tz \in \mt \cap W^u_{\ep}(\hz)$, $\ty \in W^s_\ep(\tz)$
and   a unit  vector $\xi \in E^s(\ty)$ tangent to $\mt$ at $\ty$ such that
\be
|d\omega_{\tz}(\xi_{\tz}, \eta_{\tz})| \geq \mu_0\;.
\ee

Set $\ty_1 = \ty$. Let $0 < \delta   <  \min \{\ep/2, \frac{\mu_0}{32\, C}\}$ (some additional condition on $\delta$ will be imposed later).
Since $\xi$ is tangent to $\mt$ at $\ty$, there exists $w\in E^s(\ty; \delta)$  such that $\ty_2 = \exp^s_{\ty}(w) \in \mt$ and
\be
\|w/\|w\| - \xi \| < \delta \;.
\ee
Assuming $\ep_0$ (and therefore $\ep$) is sufficiently small, there exist $w_1,w_2 \in E^s(\tz;\ep)$ such that
$\ty_i = \exp^s_{\tz}(w_i)$ for $i = 1,2$. Moreover we will assume\footnote{Using local coordinates on the Riemann manifold
$X = W^s_{\ep_0}(\tz)$ we can identify $E^s(\ty_1) = T_{\ty_1}X$ and $E^s(\tz) = T_{\tz} X$.  
Given $\omega > 0$, we can take $\ep > 0$ so small that $\|P_{y,\tz} - I\| < \omega/2$ for any $y \in X$ with $d(y,\tz) < \ep$,
where $P_{y,\tz}$ is the operator of parallel translation along the geodesic in $X$ from $y$ to $\tz$. Then,
given $\ty_1 = \exp^s_{\tz}(w_1)$ with $\|w_1\| < \ep$, a simple calculation shows that we can choose $\delta > 0$ 
so small that $\|(\exp^s_{\tz})^{-1}(\exp^s_{\ty_1}(\zeta)) - (w_1 + \zeta)\| \leq \omega \|\zeta\|$ for every $\zeta \in T_{\ty_1}X$ with
$\|\zeta\| < \delta$.}
that $C > 0$ is taken sufficiently large and $\ep > 0$ and then $\delta > 0$ sufficiently small 
so that  $\|w_2 - (w_1 + \|w\|\, \xi_{\tz})\| \leq \frac{\mu_0}{8 C}\, \|w\|$.
Fix $w$ with this property and (6.3), and set $\tde = \frac{\mu_0 \|w\|}{10}$.

Given any $z\in \mt \cap W^u_{\ep}(\tz)$, for $i = 1,2$ there exists a  unique $w_i(z) \in E^s(z)$ such that
$\pi_{\ty_i}(z) = \exp^s_z(w_i(z))$. Clearly $w_i(z)$ is a continuous function of $z$ with $w_i(\tz) = w_i$. 
Take $0 < \tep \leq \min \left\{ \ep , \frac{\mu_0}{16 C \|w\|} \right\}$
so small that the conclusion of Lemma 6.3 holds and moreover
$\|w_i(z)\| < \ep$ and $\|w_i\|/2 \leq \|w_i(z)\| \leq 2\|w_i\|$  for all $z\in \mt \cap W^u_{\tep}(\tz)$ and $i = 1,2$.

We will now use parallel translation on the Riemann manifold $W^u_{\ep_0}(z_0)$. For any
$z  \in W^u_{\tep}(\tz)$ let $\Gamma_z  : E^u(z) \longrightarrow E^u (\tz)$
be the parallel translation along the geodesic in $W^u_{\ep_0} (z_0)$ from $z$ to $\tz$. Then $\Gamma_z$
is an isometry which is Lipschitz in $z$  and $\Gamma_{\tz} = \id$. Since the form $d\omega_z$ is $C^1$ in 
$z$ and $w(z) \to w$ as $z\to \tz$, taking $\tep > 0$ sufficiently small,
for any $z  \in W^u_{\ep}(\tz) \cap \mt$ with $d(z,\tz) < \tep$ we have
\be
|d\omega_z( v , w_i(z)) - d\omega_{\tz}( \Gamma_z(v) ,w_i)| \leq \tde \, \|v\|
\quad, \quad v \in E^u(z; \tep)\:, \: i = 1,2 \;.
\ee
Moreover we can take $\tep > 0$ so small that $\|\Gamma_z(\eta_z) - \eta_{\tz} \| \leq \frac{\mu_0}{4 C} $
for any $z  \in W^u_{\tep}(\tz) \cap \mt$.

\def\teta{\tilde{\eta}}

Let $z  \in \mt \cap W^u_{\ep_0}(\tz)$, $d(z,\tz) < \tep$, and let $v\in E^u(\tz)$ , 
be such that $\exp^u_z(v) \in \mt$, $\|v\| < \tep$,  and $\la v/\|v\| , \eta_z\ra  \geq \theta_0$. 
Setting $\tv = \Gamma_z(v)$ and $\teta = \Gamma_z(\eta_z)$, 
we have $\|\tv\| = \|v\|$, $\|\teta\| = 1$, $\la \tv, \teta\ra = \la v, \eta_z\ra$ and 
$\| \teta - \eta_{\tz}\| \leq \frac{\mu_0}{4 C}$. Thus, 
$\la \tv/\|\tv\| , \teta \ra  = \la v/\|v\|, \eta_z \ra \geq \theta_0$, which combined with (6.1) gives 
$\|\tv /\|\tv\| - \teta \|^2 = 2 -2 \la \tv / \|\tv\| ,  
\teta \ra \leq 2 (1-\theta_0) \leq \frac{\mu^2_0}{64 C^2}$, so
$\|\tv/\|\tv\| - \teta\| < \frac{\mu_0}{8 C}$.
Using  Lemma 6.3, (6.4), (6.3) and $\tde < \frac{\mu_0}{10}$ (which follows from the choice of $\tde$) we now get
\begin{eqnarray*}
&    & |\Delta( \exp^u_z(v), \pi_{\ty_1}(z)) - \Delta( \exp^u_z(v), \pi_{\ty_2}(z)) |  
 =   |\Delta( \exp^u_{z}(v), \exp^s_z(w_1(z)) -  \Delta( \exp^u_{z}(v), \exp^s_z(w_2(z)) | \\
&    & \geq   |d\omega_z(v , w_1(z)) - d\omega_z(v , w_2(z)) | - \tde\, \|v\|\, (\|w_1(z)\| + \|w_2(z)\|)\\
&    & \geq    |d\omega_{\tz}( \tv , w_1) - d\omega_{\tz}( \tv , w_2)| - 2\tde\, \|v\|   
 - 2 \ep\, \tde\, \|v\| 
\geq |d\omega_{\tz}( \tv , w_1 - w_2)| - 4\tde\, \|v\|  \\
&   & \geq \|w\|\,\|v\|\,\left[ |d\omega_{\tz}( \tv/\|\tv\| , \xi_{\tz})| 
- \frac{\mu_0}{8}\, \right] - 4\tde\, \|v\|  
\geq \|w\|\,\|v\|\,\left[ |d\omega_{\tz}( \teta , \xi_{\tz})| - \frac{\mu_0}{4}\, \right] - 4\tde\, \|v\|  \\
&   & \geq \|w\|\,\|v\|\,\left[ |d\omega_{\tz}( \eta_{\tz} , \xi_{\tz})| 
- \frac{\mu_0}{2}\, \right] - 4\tde\, \|v\|  
= \left(\frac{\mu_0\, \|w\|}{2} - 4\tde\right)\, \|v\|  =  \tde \, \|v\|\;.
\end{eqnarray*}
This proves the lemma. \endofproof

\section{Regular distortion along unstable manifolds}
\setcounter{equation}{0}

In this section we briefly describe the results in \cite{kn:St3} which give sufficient conditions
for a flow over a basic set to have regular distortion along unstable manifolds. As we mentioned
in section 1, there are good reasons to believe that this condition is satisfied for
a very general class of flows on basic sets (perhaps even always).

As before, let $M$ be a $C^2$ complete Riemann manifold, $\phi_t$  be a  $C^2$  flow on $M$, and let 
$\mt$ a basic set for $\phi_t$. 

Assume that $\phi_t$ and $\mt$ satisfy the following 
{\it lower unstable pinching condition}:

\ms

\noindent
{\sc (LUPC)}:  {\it There exist  constants $C > 0$ and $0 < \alpha \leq \beta < \alpha_2 \leq \beta_2$,
and for  every $x\in \mt$ constants $\alpha_1(x) \leq \beta_1(x)$ with $\alpha \leq \alpha_1(x) 
\leq \beta_1(x) \leq \beta$
and $2\alpha_1(x) - \beta_1(x) \geq \alpha$ and a $d\phi_t$-invariant splitting 
$E^u(x) = E^u_1(x)\oplus E^u_2(x)$,
continuous with respect to $x \in \mt$, such that
$$\frac{1}{C} \, e^{\alpha_1(x) \,t}\, \|u\| \leq \| d\phi_{t}(x)\cdot u\| 
\leq C\, e^{\beta_1(x)\,t}\, \|u\|
\quad, \quad  u\in E^u_1(x) \:\:, t > 0 \;,$$
and}
$$\frac{1}{C} \, e^{\alpha_2 \,t}\, \|u\| \leq \| d\phi_{t}(x)\cdot u\| \leq C\, e^{\beta_2\,t}\, \|u\|
\quad, \quad  u\in E^u_2(x) \:\:, t > 0 \;.$$

\ms

In (LUPC) the lower part of the spectrum of $d\phi_t$ over $E^u$ is (point-wisely) pinched, however 
there is no restriction on the rest of the spectrum, except that it should be uniformly separated 
from the lower part.

Under the above condition the distribution $E^u_2(x)$ ($x\in \mt$) is integrable 
(see e.g. \cite{kn:Pes}), so (assuming $\ep_0 > 0$ is
small enough) there exists a $\phi_t$-invariant family $W^{u,2}_{\ep_0}(x)$ ($x\in \mt)$ of $C^2$ 
submanifolds of $W^u_{\ep_0}(x)$ such that $T_x(W^{u,2}_{\ep_0}(x)) = E^u_2(x)$ for all $x\in \mt$. 
Moreover (see Theorem 6.1
in \cite{kn:HPS} or the proof of Theorem B in \cite{kn:PSW}), for any $x\in \mt$, the map 
$\mt \cap W^u_{\ep_0}(x) \ni y \mapsto E^u_2(y)$ is $C^1$. However in general the distribution  
$E^u_1(x)$ ($x\in \mt$)  does not have to be integrable (see  \cite{kn:Pes}).

We now make the  additional assumption that  $E^u_1(x)$ ($x\in \mt$) is integrable:

\ms

\noindent
{\sc (I)}:  {\it There exist $\ep_0 > 0$ and a continuous $\phi_t$-invariant family 
$W^{u,1}_{\ep_0}(x)$ ($x\in \mt)$ of $C^2$ submanifolds
of $W^u_{\ep_0}(x)$ such that $T_x(W^{u,1}_{\ep_0}(x)) = E^u_1(x)$ for all $x\in \mt$, 
and moreover for any $\ep > 0$ and 
any $x\in \mt$, $\mt\cap W^u_{\ep}(x)$ is {\bf not} contained in  $W^{u,2}_{\ep}(x)$. }

\ms

Roughly speaking, the latter means that the distribution $E^u_1(x)$ is significantly involved
in the dynamics of the flow over $\mt$. 

The main result in \cite{kn:St3} is the following.

\bs

\noindent
{\bf Theorem 7.1.} {\it Let $\phi_t$ and $\mt$ satisfy the conditions {\rm (LUPC)} 
and {\rm (I)}. Then $\phi_t$
has a regular distortion along unstable manifolds over $\mt$.}

\bs

A simplified case is presented by the following {\it pinching condition}:

\ms

\noindent
{\sc (P)}:  {\it There exist  constants $C > 0$ and $\beta \geq \alpha > 0$ such that for every $x\in \mt$ we have
$$\frac{1}{C} \, e^{\alpha_x \,t}\, \|u\| \leq \| d\phi_{t}(x)\cdot u\| \leq C\, e^{\beta_x\,t}\, \|u\|
\quad, \quad  u\in E^u(x) \:\:, t > 0 \;,$$
for some constants $\alpha_x, \beta_x > 0$ depending on $x$ but independent of $u$ and $t$ with
$\alpha \leq \alpha_x \leq \beta_x \leq \beta$ and $2\alpha_x - \beta_x \geq \alpha$ for all $x\in \mt$.}

\ms

Clearly the condition (P) is (LUPC) in the special case when $E^u_2(x) = 0$ for all $x\in \mt$.  
Notice that when the local unstable manifolds are one-dimensional the condition (P) is always satisfied. 
In higher dimensions a well-known example when (P) holds is the geodesic flow on a manifold with 
strictly negative sectional curvature satisfying the so called $\frac{1}{4}$-pinching condition 
(see \cite{kn:HP}). For open billiard flows (in any dimension) it was shown in \cite{kn:St2} that 
if the distance between the scatterers is large compared with the maximal sectional curvature of 
the boundaries, then the condition (P) is satisfied over the non-wandering set.

As a special case of Theorem 7.1, it is shown in \cite{kn:St3} that if
$\phi_t$ satisfies the condition (P) on $\mt$, then $\phi_t$
has a regular distortion along unstable manifolds over $\mt$.


\section{Appendix: Proof  of Lemma 5.4}
\setcounter{equation}{0}

 (a) Let  $u, u' \in \hU_i $ for some $i = 1, \ldots,k$ and let $m \geq 1$ be an integer.  
Given $v \in \hU$ with $\sigma^m(v) = u$, let $C[\ii] = C[i_0, \ldots,i_m]$ be the 
cylinder of length $m$ containing $v$ (see the beginning of section 3).   Since the sequence 
$\ii = [i_0, \ldots,i_m]$ is admissible, the Markov property implies $i_m = i$ and $\sigma^m(\hC[\ii]) = \hU_i$.
Moreover, $\sigma^m : \hC[\ii] \longrightarrow \hU_i $ is a homeomorphism, so
there  exists a unique $v' = v'(v)\in \hC[\ii]$ such that $\sigma^m(v') = u'$. By (2.1),
$d (\sigma^j(v),\sigma^j(v'(v))) \leq \frac{1}{c_0\, \gamma^{m-j}}\, d (u,u')$ for all $j = 0,1, \ldots, m-1\;.$
This and (2.2) imply
$$|\fa_m(v) - \fa_m(v')| \leq
\sum_{j=0}^{m-1} \frac{\Lip(\fa)}{c_0\, \gamma^{m-j}} \, d (u,u')
 \leq  \frac{ T}{c_0\, (\gamma-1)}\, d (u,u')
\leq \frac{ T}{c_0\, (\gamma-1)}\, D (u,u')\;.$$

Also notice that if $D(u,u') = \diam(\cc')$ for some cylinder $\cc' = C[i_m, i_{m+1}, \ldots, i_p]$,
then $v,v'(v) \in \cc'' = C[i_0,i_1, \ldots,i_p]$ for some cylinder $\cc''$ with
$\sigma^m(\cc'') = \cc'$, so 
$D(v,v'(v)) \leq \diam(\cc'')  \leq \frac{1}{c_0\, \gamma^m} \, \diam(\cc') 
= \frac{D(u,u')}{c_0\, \gamma^m} \;.$

Using the above, $\diam (U_i) \leq 1$, the definition of $\ma$, and the fact that 
$\ma 1 = 1$ (hence $\ma^m 1 = 1$), and assuming $A_0 \geq e^{\frac{T}{c_0\, (\gamma-1)}}/c_0$, we get
\begin{eqnarray*}
&        &\frac{|(\ma^m H)(u) - (\ma^m H)(u')|}{\ma^m H (u')} 
 =      \frac{\di \left| \sum_{\sigma^m v = u} e^{\fa_m(v)}\, H(v) -  
\sum_{\sigma^m v = u} e^{\fa_m(v'(v))}\, H(v'(v)) \right|}{\ma^m H (u')} \\
& \leq & \frac{\di \left| \sum_{\sigma^m v = u} e^{\fa_m(v)}\, (H(v) -  H(v'(v)))\right|}{\ma^m H (u')}   +
\frac{\di \sum_{\sigma^m v = u}  \left|e^{\fa_m(v)}-  e^{\fa_m(v'(v))}\right| \, H(v'(v))}{\ma^m H (u')} \\
& \leq & \frac{\di  \sum_{\sigma^m v = u} e^{\fa_m (v)}\, B\, H(v'(v))\, D (v,v'(v))}{\ma^m H (u')} 
 + \frac{\di \sum_{\sigma^m v = u} \left| e^{\fa_m(v)- \fa_m(v'(v))} -1 \right|  \,
e^{\fa_m(v'(v))}\, H(v'(v))}{\ma^m H (u')}\\ 
& \leq &  e^{\frac{T}{c_0\, (\gamma-1)}}\, \frac{B\, D (u,u')}{c_0\gamma^m}
 + e^{\frac{T}{c_0\, (\gamma-1)}}\, \frac{T}{c_0\, (\gamma-1)}\, D (u, u')
 \leq  A_0 \, \left[ \frac{B}{\gamma^m} + \frac{T}{\gamma-1}\right]\, D (u,u')\;.
\end{eqnarray*}

\ms

(b) The proof of this part is very similar to the above and we omit it. \endofproof

\ms


\footnotesize

\noindent
{\bf Acknowledgements.} Thanks are due to Keith Burns, Boris Hasselblatt, Charles Pugh and 
Amie Wilkinson for various kind of information they have provided to me. The work on this article 
(together with \cite{kn:St3}) started a long time ago  during my visit to  the Erwin Schr\"odinger 
Institute in Vienna for the Program on Scattering Theory in 2001
organized by Vesselin Petkov, Andras Vasy and Maciej Zworski. Just by chance a significant part 
of the work on the latest revision of the paper was done during my visit to ESI for the Program 
on Hyperbolic Dynamical Systems in 2008
organized by Domokos Sz\'asz, Lai-Sang Young and Harald Posch.
Thanks are due to the organizers of these events and to the staff of ESI for their hospitality and support.


\bs

{\sc University of Western Australia, Crawley WA 6009, Australia}

{\sc\it E-mail address:} stoyanov@maths.uwa.edu.au

\end{document}